\documentclass[10pt,a4paper,final,english]{amsart}
\usepackage[utf8]{inputenc}
\usepackage{babel}
\usepackage{csquotes}
\usepackage[all]{xy}
\usepackage{graphicx}
\usepackage{color}
\usepackage{xspace}
\usepackage{amsmath}
\usepackage{amsthm}
\usepackage{amssymb}
\usepackage{mathtools}
\usepackage{tensor} 
\usepackage{scalerel,stackengine}
\usepackage{doi}
\usepackage[style=alphabetic,
            sorting=nyt,
            url=false,
            isbn=false,
            doi=false,
            clearlang=true,
            maxnames=4,
            giveninits=true,
            backend=biber]{biblatex}
\hypersetup{hidelinks}

\newcommand{\baseRing}[1]{\ensuremath{\mathbb{#1}}}
\newcommand{\R}{\baseRing{R}}
\newcommand{\FP}[2]{\tensor[_{{#1}}]{\times}{_{{#2}}}} 
\newcommand{\gpoidarrows}{\ensuremath{\rightrightarrows}}

\newcommand{\FBS}{\ensuremath{{Fbs}}\xspace}
\newcommand{\setC}{\ensuremath{Set}\xspace}
\newcommand{\lLgpdC}{\ensuremath{lLGpd}\xspace}
\newcommand{\LalgbdC}{\ensuremath{LAlgbd}\xspace}
\newcommand{\VBC}{\ensuremath{Vb}\xspace}
\newcommand{\VFunctor}{\ensuremath{V}}
\newcommand{\LieFunctor}{\ensuremath{Lie}}

\DeclareMathOperator{\Ad}{Ad}
\DeclarePairedDelimiter{\norm}{\lVert}{\rVert}

\stackMath
\newcommand\widecheck[1]{%
  \savestack{\tmpbox}{\stretchto{%
      \scaleto{%
        \scalerel*[\widthof{\ensuremath{#1}}]{\kern-.6pt\bigwedge\kern-.6pt}%
        {\rule[-\textheight/2]{1ex}{\textheight}}
      }{\textheight}%
    }{0.5ex}}%
  \stackon[1pt]{#1}{\scalebox{-1}{\tmpbox}}%
}
\newcommand{\DC}{\ensuremath{{\mathcal{A}_d}}}
\newcommand{\DCp}[1]{\ensuremath{{\mathcal{A}_d^{#1}}}}
\newcommand{\SG}{\ensuremath{G}\xspace}
\newcommand{\jgg}{\ensuremath{\mathfrak{g}}\xspace}
\newcommand{\CC}{\ensuremath{{\mathcal{A}}}}
\newcommand{\jdef}[1]{\index{#1}\emph{#1}}
\newcommand{\HLd}[1]{\ensuremath{{h_d^{{#1}}}}}
\newcommand{\HLds}[1]{\ensuremath{\overline{h_d^{{#1}}}}}

\newcommand{\HLc}[1]{\ensuremath{{h^{{#1}}}}}
\newcommand{\jcU}{\ensuremath{\mathcal{U}}\xspace}
\newcommand{\jcV}{\ensuremath{\mathcal{V}}\xspace}
\newcommand{\jcW}{\ensuremath{\mathcal{W}}\xspace}
\newcommand{\jcZ}{\ensuremath{\mathcal{Z}}\xspace}
\newcommand{\jcG}{\ensuremath{\mathcal{G}}\xspace}
\newcommand{\jcH}{\ensuremath{\mathcal{H}}\xspace}
\newcommand{\ti}[1]{\widetilde{#1}}
\newcommand{\stext}[1]{\ensuremath{\quad\text{#1}\quad}}
\newcommand{\CB}{\ensuremath{{\mathcal{B}}\xspace}}
\newcommand{\DB}{\ensuremath{{\mathcal{B}_d}}}
\newcommand{\lie}[1]{\ensuremath{Lie(#1)}}
\newcommand{\VF}{\mathfrak{X}}
\newcommand{\VC}{\ensuremath{\mathcal{V}}\xspace}

\providecommand{\ip}[2]{\ensuremath{\left\langle{#1},{#2}\right\rangle}}
\newcommand{\del}{\ensuremath{\partial}}
\newcommand{\CE}{\ensuremath{\coloneqq}\xspace}
\newcommand{\difference}{\alpha} \newcommand{\RR}{\mathcal R}
\newcommand{\Hdr}{H_{\mathrm {dR}}} \newcommand{\openset}{\mathcal{V}}
\newcommand{\fint}{H} \newcommand{\vbundle}{\mathcal{V}}
\newcommand{\Epsilon}{\mathcal{E}}

\theoremstyle{plain}
\newtheorem{theorem}{Theorem}[section]
\newtheorem{corollary}[theorem]{Corollary}
\newtheorem{prop}[theorem]{Proposition}
\newtheorem{proposition}[theorem]{Proposition}
\newtheorem{lemma}[theorem]{Lemma}

\theoremstyle{definition}
\newtheorem{definition}[theorem]{Definition}
\newtheorem{remark}[theorem]{Remark}
\newtheorem{example}[theorem]{Example}

\newtheorem*{problem*}{Problem}
\numberwithin{equation}{section}

\begin{document}
  
\title[The integration problem]{The integration problem\\for principal connections}

\author{Javier Fern\'andez}
\address{Departamento de Matem\'atica, Instituto Balseiro, Universidad
  Nacional de Cuyo -- C.N.E.A.\\Av. Bustillo 9500, Bariloche, 8400, R{\'\i}o
  Negro, Rep\'ublica Argentina}
\email{jfernand@ib.edu.ar}

\author{Francisco Kordon}
\address{Universidad Nacional de San Mart{\'\i}n (UNSAM), Campus
  Miguelete, 25 de Mayo y Francia, CP 1650, San Mart{\'\i}n, Argentina}
\email{fkordon@unsam.edu.ar}





\begin{abstract}
  In this paper we introduce the Integration Problem for principal
  connections. Just as a principal connection on a principal $\SG$-bundle
  $\phi:Q\rightarrow M$ may be used to split $TQ$ into horizontal and
  vertical subbundles, a discrete connection may be used to split
  $Q\times Q$ into horizontal and vertical submanifolds. All discrete
  connections induce a connection on the same principal bundle via a
  process known as the Lie or derivative functor. The Integration
  Problem consists of describing, for a principal connection $\CC$,
  the set of all discrete connections whose associated connection is
  $\CC$.  Our first result is that for \emph{flat} principal
  connections, the Integration Problem has a unique solution among the
  \emph{flat} discrete connections. More broadly, under a fairly mild
  condition on the structure group $\SG$ of the principal bundle
  $\phi$, we prove that the existence part of the Integration Problem
  has a solution that needs not be unique. Last, we see that, when
  $\SG$ is abelian, given compatible continuous and discrete
  curvatures the Integration Problem has a unique solution constrained
  by those curvatures.
\end{abstract}

\thanks{This research was partially supported by grants from the
  Universidad Nacional de Cuyo (grants 06/C567 and 06/C496) and
  CONICET}

\thanks{The reference for the published version of this preprint is
  \doi{10.1016/j.geomphys.2025.105566}.}

\maketitle
  

\section{Introduction}
\label{sec:introduction}


Connections on bundles are very useful tools both in the applications
of bundles to different geometric problems ---including some areas of
Physics, such as Classical
Mechanics~\cite{bo:cendra_marsden_ratiu-lagrangian_reduction_by_stages},
Gauge
Theory~\cite{ar:daniel_viallet-the_geometrical_setting_of_gauge_theories_of_the_yang_milss_type},
and Topological
Insulators~\cite{ar:fruchart_carpentier-an_introduction_to_topological_insulators}---
and, also, in the study of the geometry of the bundle itself. Let
$\phi:Q\rightarrow M$ be a principal $\SG$-bundle. Roughly speaking, a
principal connection is a way of providing the tangent spaces $T_qQ$
with subspaces that are complementary to (the tangent space to) the
orbit of $\SG$ through $q$ in $Q$. Because of the applicability of
this notion, a principal connection has been characterized in several
different ways, as we review later.

A discrete connection on $\phi$ is also a way of providing a
complement to the $\SG$-orbits, but this time in $Q\times Q$
---instead of $TQ$--- in some appropriate sense. Again, discrete
connections can be characterized in different ways, to be recalled
later in Section~\ref{sec:connections_and_discrete_connections}.

Principal connections on a principal $\SG$-bundle
$\phi:Q\rightarrow M$ appear naturally in the study of symmetries of
some continuous-time dynamical systems defined on $Q$, where some of
the standard constructions require splitting the velocity vectors
$v_q\in T_qQ$ into horizontal and vertical parts (see, for
example,~\cite{bo:cendra_marsden_ratiu-lagrangian_reduction_by_stages}). Discrete
connections were introduced in this context as a tool used to perform
the same kind of construction for some discrete-time dynamical
systems, where the role of $TQ$ is played by $Q\times Q$ ---in a naive
sense, velocities are replaced by pairs of nearby points
(see~\cite{ar:leok_marsden_weinstein-a_discrete_theory_of_connections_on_principal_bundles},~\cite{ar:marrero_martin_martinez-discrete_lagrangian_and_hamiltonian_mechanics_on_lie_groupoids},
and~\cite{ar:fernandez_zuccalli-a_geometric_approach_to_discrete_connections_on_principal_bundles}). Interestingly,
discrete connections have also been considered in the context of
Synthetic Differential Geometry as a kind of combinatorial analogue
(on a first neighborhood of the diagonal of a manifold) of principal
connections (see
\cite{ar:kock-principal_bundles_groupoids_and_connections}, for
instance). Some other notions inspired by that of a principal
connection have been considered, for example,
in~\cite{ar:fernandes_marcut-multiplicative_ehresmann_connections}
and~\cite{laurentGengoux_stienon_xu-non_abelian_gerbes-2009}; see
Remark~\ref{rem:MECs_and_DCs} for a comparison.

It has always been known (Section 5.2
in~\cite{ar:leok_marsden_weinstein-a_discrete_theory_of_connections_on_principal_bundles})
that there is a ``derivation'' process that can be used in a principal
bundle $\phi$ to produce a principal connection $\CC$ out of a
discrete connection $\DC$, both on $\phi$. Loosely speaking, given a
discrete connection $\DC$, the associated principal connection $\CC$
is ``its derivative'', while any discrete connection $\DC$ whose
derivative is $\CC$ is called ``an integral'' of $\CC$. Besides this
construction, not much is known about the relationship between the two
types of connection, other than the fact that both have very similar
properties.

The purpose of this paper is to introduce and study aspects of the
following questions, that we collectively call \jdef{the Integration
  Problem}.

\begin{problem*}\label{pr:general_problems}
  What can we say about the ``derivation'' mapping from discrete
  connections to continuous connections on a principal bundle? More
  specifically:
  \begin{enumerate}
  \item \label{it:general_problems-surjctive} Is it surjective? That
    is, given a connection $\CC$ on a principal bundle, is there any
    discrete connection $\DC$ on that bundle that maps into $\CC$?
  \item \label{it:general_problems-injective} Is it injective? Maybe better,
    what is the structure of the set of discrete connections $\DC$
    that map to the same continuous connection $\CC$?
  \item \label{it:general_problems-make_injective} If this mapping is
    not injective, are there any additional constraints that could be
    imposed in order to make it so?
  \end{enumerate}
\end{problem*}

In Section~\ref{sec:connections_and_discrete_connections} we review
three equivalent characterizations of connections and of discrete
connections in principal bundles: using forms, horizontal liftings and
splittings of certain sequences. We also study how a derivation
functor relates the discrete to the continuous notions.

We focus our attention in Section~\ref{sec:the_flat_case} on a type of
connection that has specially nice properties: the flat connections,
both continuous and discrete. It is well known that flat principal
connections on $\phi:Q\rightarrow M$ are equivalent to morphisms of
Lie algebroids $TM\rightarrow TQ/\SG$
(see~\cite{bo:mackenzie-lie_groupoids_and_algebroids_in_differential_geometry}). It
has been recently seen that flat discrete connections are equivalent
to morphisms of local Lie groupoids
$M\times M\rightarrow (Q\times Q)/\SG$
(see~\cite{ar:fernandez_juchani_zuccalli-discrete_connections_on_principal_bundles_the_discrete_atiyah_sequence}). Then,
the Integration Problem becomes an instance of Lie Theory for Lie
groupoids and algebroids, especially Lie's Second Theorem. Using a
result
of~\cite{ar:cabrera_marcut_salazar-on_local_integration_of_lie_brackets},
we easily prove in Theorem~\ref{th:solution_of_IP-flat_case} that the
Integration Problem has an essentially unique solution, meaning that,
given a flat principal connection $\CC$, there is an essentially
unique flat discrete connection $\DC$ (unique as a germ) that
integrates $\CC$.

Leaving aside the context of flat connections, in
Section~\ref{sec:existence_of_solution_in_the_general_case} we prove
(Corollary~\ref{cor:DC_from_retraction_induces_CC}) that arbitrary
principal connections can be integrated provided that the principal
bundle is equipped with a smooth retraction (see Chapter 4
of~\cite{bo:absil_mahony_sepulchre-optimization_algorithms_on_matrix_manifolds}). Then
we prove that such retractions do exist if, for instance, the
structure group of the principal bundle is the Cartesian product of a
compact Lie group and a vector space
(Corollary~\ref{cor:existence_of_integrals_for_compact_G}).  Thus, we
have established that, under appropriate conditions, the existence
part of the Integration Problem has a solution. Moreover, as we show
in an example, this solution needs not be unique.

Last, in Section~\ref{sec:the_abelian_structure_group_case}, we
consider the special case where the structure group $\SG$ of the
principal bundle is abelian. In this setting we see that the
Integration Problem for arbitrary curvatures has a solution if and
only if the curvatures of the principal and the discrete connections
are related in a specific way. Then we prove that, when that relation
occurs, the Integration Problem for connections whose curvatures
satisfy that relationship admits an essentially unique solution. This
result is compatible with the result from
Section~\ref{sec:the_flat_case} where the two curvatures are chosen to
be flat, which, indeed, satisfies the said relation. Therefore, we
proved the existence and uniqueness of the solution to the Integration
Problem when the structure group of the bundle is abelian.

The general Integration Problem remains open. The partial solutions
that we provide in this paper might serve as indications that this is
a rich and nuanced topic, worthy of further study.

Finally, we wish to thank A. Cabrera and I. M\u{a}rcu\c{t} for
discussing with us some of their results
in~\cite{ar:cabrera_marcut_salazar-on_local_integration_of_lie_brackets}
as well as for their interest in our work.

\subsubsection*{Notation}
\label{sec:notation}

Throughout the paper, if the Lie group $G$ acts on the manifold $X$ on
the left we denote the action by $l^X$ and the lifted action on $TX$
by $l^{TX}$, so that, if $g\in G$ and $\delta x\in T_xX$, we have
$l^{TX}_g(\delta x) \CE  T_xl_g^Q(\delta x)$. In this context,
$\pi^{X,G}:X\rightarrow X/G$ is the quotient map and
$\tau_X:TX\rightarrow X$ is the canonical projection.


\section{Connections and discrete connections}
\label{sec:connections_and_discrete_connections}

In this section we recall some of the characterizations and basic
properties of connections and discrete connections on a principal
bundle.


\subsection{Connections on principal bundles}
\label{sec:CC_on_principal_bundles}

As we mentioned in the Introduction, connections on principal bundles
may be described in several equivalent ways. Next we review three such
approaches that will suffice for the purpose of this paper. For more
details, see~\cite{bo:kobayashi_nomizu-foundations-v1}
and~\cite{bo:mackenzie-lie_groupoids_and_algebroids_in_differential_geometry}.

Let $\phi:Q\rightarrow M$ be a smooth principal (left) $\SG$-bundle
and $\jgg\CE \lie{G}$; we denote the corresponding $\SG$-action on $Q$
by $l^Q$. A connection on $\phi$ is determined by a $\jgg$-valued
$1$-form $\CC$ on $TQ$ satisfying $\CC(\xi_Q(q)) = \xi$ for all
$q\in Q$ and $\xi\in\jgg$, where
\begin{equation*}
  \xi_Q(q) \CE  \frac{d}{dt}\bigg|_{t=0} l^Q_{\exp(t\xi)}(q),
\end{equation*}
and $\CC(l^{TQ}_g(v)) = \Ad_g(\CC(v))$ for all $v\in TQ$ and
$g\in\SG$, where $\Ad$ is the adjoint action of $\SG$ on $\jgg$. Such
$\CC$ is called a \jdef{connection $1$-form}. The set of all
connection $1$-forms on $\phi$ will be denoted by $\Sigma_C$.

\begin{example}\label{ex:trivial_bundle:cc}
  Given a Lie group $G$ and a manifold $M$ we define
  $Q\CE M \times G$. The induced left $G$-action on $Q$, defined
  by $l_g^Q(m,g')\CE (m,gg')$ for $m\in M$ and $g,g'\in G $, makes the
  projection onto the first component $p_1:Q\rightarrow M$ a principal
  $G$-bundle.

  Given a connection form $\CC$ on $p_1$ we may define the $1$-form
  $\omega\in\Omega^1(M,\jgg)$, which is often called the \jdef{local
    expression} of~$\CC$, by precomposition with the differential of
  the trivial section of $p_1$:
  \begin{equation*}
    \omega_m(\delta m) 
    \CE \CC_{(m,e)}(\delta m \oplus 0)
    \quad \text{for } m\in M \text{ and } \delta m\in T_mM.
  \end{equation*}
  Conversely, each form $\omega\in\Omega^1(M,\jgg)$ gives a connection
  $\CC$ on $p_1$ defined by
  $\CC(\delta m + \delta g) \CE \omega(\delta m)
  +Tr^G_{g^{-1}}(\delta g)$ for $\delta m\in T_mM$ and
  $\delta g\in T_gG$, where $r^G_h$ denotes the right multiplication
  on $G$ by $h\in G$. These two constructions are reciprocal inverses.
\end{example}

Alternatively, a connection on $\phi$ may be described using a smooth
morphism of vector bundles $\HLc{}:\phi^*TM\rightarrow TQ$ that is
$\SG$-equivariant for the $\SG$-actions $l^{TQ}$ and
$l^{Q\times TM}_g(q,\delta m)\CE (l^Q_g(q),\delta m)$, and
$(\tau_Q,T\phi)\circ \HLc{} = id_{\phi^*TM}$. Such $\HLc{}$ is called a
\jdef{horizontal lift}. The set of all horizontal lifts $\HLc{}$ on
$\phi$ is denoted by $\Sigma_H$.

Given $\CC\in\Sigma_C$, we define $\HLc{}:\phi^*TM\rightarrow TQ$ by
\begin{equation}\label{eq:F_{CH}-def}
  \HLc{}(q,\delta m) \CE  \delta q - (\CC(\delta q))_Q(q),
\end{equation}
where $\delta q\in T_qQ$ is such that $T_q\phi(\delta q)=\delta
m$. The following result is fairly straightforward to check.
\begin{prop}
  The map $F_{CH}:\Sigma_C\rightarrow \Sigma_H$ where
  $F_{CH}(\CC)\CE \HLc{}$ with $\HLc{}$ given
  by~\eqref{eq:F_{CH}-def} is well defined. Furthermore, $F_{CH}$
  is a bijection.
\end{prop}

Yet another description of connections on $\phi$ may be obtained by
considering the splittings of the following exact sequence in the
category of vector bundles over $M$, $\VBC_M$, the subcategory of the
category of vector bundles $\VBC$ whose objects have base $M$ and
morphisms over $id_M$. Consider the sequence
\begin{equation}\label{eq:AS-def}
  \xymatrix{ {\{0\}} \ar[r] & {\ti{\jgg}} \ar[r]^-{\phi_1} &
    {(TQ)/\SG} \ar[r]^-{\phi_2} & {TM} \ar[r] & {\{0\}}
  },
\end{equation}
where $\ti{\jgg}$ is the \jdef{adjoint bundle} $(Q\times \jgg)/\SG$
with $\SG$ acting by $l^Q$ and by $\Ad$; $\phi_1([q,\xi])\CE [\xi_Q(q)]$
and $\phi_2([v_q])\CE T_q\phi(v_q)$. The sequence~\eqref{eq:AS-def} is
exact (Proposition 3.2.3
of~\cite{bo:mackenzie-general_theory_of_lie_groupoids_and_algebroids})
and is known as the \jdef{Atiyah Sequence} associated to the principal
bundle $\phi$. It is well known that right splittings
of~\eqref{eq:AS-def} also determine principal connections on
$\phi$. The set of all right splittings of~\eqref{eq:AS-def} is
denoted by $\Sigma_R$.

Given $\HLc{}\in \Sigma_H$, we define $s:TM\rightarrow (TQ)/\SG$ by
\begin{equation}\label{eq:s_from_HL-def}
  s(\delta m) \CE  [\HLc{}(q,\delta m)] 
\end{equation}
for any $\delta m\in T_{\phi(q)}M$. Then, we have the following result.
\begin{prop}
  The map $F_{HR}:\Sigma_H\rightarrow \Sigma_R$ where
  $F_{HR}(\HLc{})\CE s$ with $s$ given
  by~\eqref{eq:s_from_HL-def} is well defined. Furthermore,
  $F_{HR}$ is a bijection.
\end{prop}

We define the map $F_{CR}:\Sigma_C\rightarrow \Sigma_R$ as
$F_{CR} \CE  F_{HR}\circ F_{CH}$. Naturally, $F_{CR}$ is a
bijection. All together, we have the following commutative diagram in
the category of sets
\begin{equation*}
  \xymatrix{
    {} & {\Sigma_R}  & {}\\
    {\Sigma_C} \ar[rr]_{F_{CH}}^{\sim} \ar[ur]^{F_{CR}}_{\sim}& {} &
    {\Sigma_H} \ar[ul]_{F_{HR}}^{\sim} 
  }
\end{equation*}


\subsection{Discrete connections on principal bundles}
\label{sec:DC_on_principal_bundles}

Roughly speaking, discrete connections on the principal $\SG$-bundle
$\phi:Q\rightarrow M$ are similar to connections but substituting $TQ$
by $Q\times Q$ and $\jgg$ by $\SG$. They also require an additional piece of
information: while connections are globally defined, discrete
connections may only exist semi-locally, as we see next. For more
details,
see~\cite{ar:fernandez_zuccalli-a_geometric_approach_to_discrete_connections_on_principal_bundles}
and~\cite{ar:fernandez_juchani_zuccalli-discrete_connections_on_principal_bundles_the_discrete_atiyah_sequence}.

An open subset $\jcU\subset Q\times Q$ is said to be of
\jdef{$D$-type} if it is $\SG\times \SG$-invariant and contains the
diagonal $\Delta_Q$.  If $\jcU$ is of $D$-type, it is said to be
\jdef{symmetric} if $(q,q')\in\jcU\iff (q',q)\in \jcU$. When $\jcU$ is
of $D$-type, we define
$\jcU''\CE (\phi\times \phi)(\jcU)\subset M\times M$ and
$\jcU'\CE (id_Q\times \phi)({\jcU})\subset Q\times M$. Both $\jcU'$ and
$\jcU''$ are open subsets.

One way of characterizing a discrete connection on $\phi$ is via a
discrete connection form, as follows. Let $\jcU$ be of $D$-type; a
smooth function $\DC:\jcU\rightarrow \SG$ is a \jdef{discrete
  connection form} if $\DC(q,q)=e$ for all $q\in Q$ and
$\DC(l^Q_{g}(q),l^Q_{g'}(q')) = g'\DC(q,q')g^{-1}$ for all
$(q,q')\in\jcU$ and $g,g'\in \SG$. The set of all discrete connection
forms on $\phi$ with domain $\jcU$ is denoted by $\Sigma_C^d(\jcU)$.

\begin{example}\label{ex:trivial_budle:dc}
  Let $Q=M\times G$ and $p_1:Q\rightarrow M$ be the trivial bundle of
  Example~\ref{ex:trivial_bundle:cc}. Let $\jcU''\subset M\times M$ be
  an open subset containing the diagonal $\Delta_M$ and $C:\jcU''\rightarrow G$
  be a smooth function such that $C(m,m)=e$ for all $m\in M$. The
  $G$-valued function $\DC$ with domain $(p_1\times
  p_1)^{-1}(\jcU'')$ defined by
  \begin{equation}\label{eq:trivial_bundle:dc}
    \DC\left( (m_0,g_0), (m_1,g_1) \right)
    \CE g_1 C(m_0,m_1)g_0^{-1}
  \end{equation}
  is a discrete connection form on $p_1$, as one can readily check. 
  Conversely, any connection form $\DC$ defined on a $D$-type open subset 
  $\jcU$ of $Q\times Q$ gives rise to a smooth function $C:(p_1\times
  p_1)(\jcU)\rightarrow G$ such that $C(m,m)=e$ for any $m\in M$ and
  \eqref{eq:trivial_bundle:dc} is valid.
\end{example}

An alternative way of characterizing a discrete connection on $\phi$
with domain $\jcU$ is using a discrete horizontal lift. A smooth map
$\HLd{}:\jcU'\rightarrow Q\times Q$ is a \jdef{discrete horizontal
  lift} if $\HLd{}(l^Q_g(q),m) = l^{Q\times Q}_g(\HLd{}(q,m))$ for all
$(q,m)\in\jcU'$ and $g\in \SG$,
$(id_Q\times \phi) \circ \HLd{} = id_{\jcU'}$, and
$\HLd{}(q,\phi(q)) = (q,q)$ for all $q\in Q$. The set of all discrete
horizontal lifts on $\phi$ with domain $\jcU'$ is denoted by
$\Sigma_H^d(\jcU)$. It is convenient to consider the second component
of a discrete horizontal lift: $\HLds{}:Q\times M\rightarrow Q$ given
by $\HLds{}\CE p_2\circ \HLd{}$.

Given $\DC\in \Sigma_C^d(\jcU)$ for $\jcU$ of $D$-type, we define
$\HLd{}:\jcU'\rightarrow Q\times Q$ by
\begin{equation}\label{eq:F_{CH}^d-def}
  \HLd{}(q,m)\CE (q,l^Q_{\DC(q,q')^{-1}}(q')) \stext{ for any } q'\in \phi^{-1}(m).
\end{equation}

\begin{prop}
  The map $F_{CH}^d:\Sigma_C^d(\jcU)\rightarrow \Sigma_H^d(\jcU)$ where
  $F_{CH}^d(\DC)\CE \HLd{}$ with $\HLd{}$ given
  by~\eqref{eq:F_{CH}^d-def} is well defined. Furthermore, $F_{CH}^d$
  is a bijection.
\end{prop}

\begin{proof}
  It can be derived from Theorems 3.6 and 4.6
  of~\cite{ar:fernandez_zuccalli-a_geometric_approach_to_discrete_connections_on_principal_bundles}.
\end{proof}

Just as in the continuous case, yet another description of discrete
connections on $\phi$ is possible, in terms of right splittings of a
certain sequence. Instead of working in the $\VBC$ category as in the
continuous case, here we consider the category $\FBS$ whose objects
are smooth fiber bundles with a given global section and whose
morphisms are smooth maps that commute with the bundle projections and
map the given sections one to the other. $\FBS_M$ is the subcategory
of $\FBS$ whose objects are fiber bundles over the manifold $M$ and
whose morphisms are smooth maps over $id_M$.

Consider the sequence
\begin{equation}\label{eq:DAS-def}
  \xymatrix{
    {\ti{\SG}} \ar[r]^-{F_1} & {(Q\times Q)/\SG} \ar[r]^-{F_2} & {M\times M}
  }
\end{equation}
where $\ti{\SG}$ is the \jdef{conjugate bundle} $(Q\times\SG)/\SG$ for
$\SG$ acting on $Q$ by $l^Q$ and on itself by conjugation. All three
spaces can be seen as fiber bundles over $M$ via the maps induced by
the projection onto the first component. Each fiber bundle is equipped
with the global section:
\begin{equation}\label{eq:DAS_spaces-sigma-def}
  \sigma_{\ti{\SG}}(m)\CE [(q,e)],\quad \sigma_{(Q\times Q)/\SG}(m)\CE [(q,q)],\quad
  \sigma_{M\times M}(m)\CE (m,m),
\end{equation}
where $q\in \phi^{-1}(m)$ is arbitrary. The maps $F_1$ and $F_2$ are
defined by
\begin{equation}\label{eq:DAS_maps-def}
  F_1([q,g]) \CE  [(q,l^Q_g(q))] \stext{ and } F_2([q,q']) \CE  (\phi(q),\phi(q')).
\end{equation}
The sequence~\eqref{eq:DAS-def} is known as the \jdef{Discrete Atiyah
  Sequence}. It is easy to check that it is a sequence in the $\FBS_M$
category.

Let $(E,\sigma),(E',\sigma') \in \FBS_M$ and $\jcV\subset E'$
be an open subset containing $\sigma'(M)$; a map
$H:\jcV\rightarrow E$ is a \jdef{semi-local morphism} if
$\phi\circ H = \phi'|_{\jcV}$ and $\sigma = H\circ \sigma'$. Let
$(E'',\sigma'')\xrightarrow{f_1} (E,\sigma) \xrightarrow{f_2}
(E',\sigma')$ be a sequence in the $\FBS_M$ category. A semi-local
morphism $H:\jcV\rightarrow E$ is a \jdef{semi-local right
  splitting} of the sequence if $f_2\circ H = id_{\jcV}$. If $\jcU$ is
of $D$-type, the set of all semi-local right splittings of the
Discrete Atiyah Sequence with domain $\jcU''$ is denoted by
$\Sigma_R^d(\jcU)$.

Let $\jcU$ be of $D$-type and $\HLd{}\in\Sigma_H^d(\jcU)$. We define
$s_d:\jcU''\rightarrow (Q\times Q)/\SG$ by
\begin{equation}\label{eq:s_d_from_HLd-def}
  s_d(m,m')\CE [\HLd{}(q,m')] \stext{ for any } q\in \phi^{-1}(m).
\end{equation}

\begin{prop}
  The map $F_{HR}^d:\Sigma_H^d(\jcU)\rightarrow \Sigma_R^d(\jcU)$
  where $F_{HR}^d(\HLd{})\CE s_d$ with $s_d$ given
  by~\eqref{eq:s_d_from_HLd-def} is well defined. Furthermore,
  $F_{HR}^d$ is a bijection.
\end{prop}

\begin{proof}
  That $F_{HR}^d$ is well defined is a matter of routine checking. The
  second assertion is Proposition 3.15
  in~\cite{ar:fernandez_juchani_zuccalli-discrete_connections_on_principal_bundles_the_discrete_atiyah_sequence}.
\end{proof}

Last, if $\jcU$ is of $D$-type, we define
$F_{CR}^d:\Sigma_C^d(\jcU)\rightarrow \Sigma_R^d(\jcU)$ as
$F_{CR}^d \CE  F_{HR}^d\circ F_{CH}^d$ and, naturally, $F_{CR}^d$ is a
bijection. We have the following commutative diagram in the category
of sets
\begin{equation*}
  \xymatrix{
    {} & {\Sigma_R^d(\jcU)}  & {}\\
    {\Sigma_C^d(\jcU)} \ar[rr]_{F_{CH}^d}^{\sim} \ar[ur]^{F_{CR}^d}_{\sim}& {} &
    {\Sigma_H^d(\jcU)} \ar[ul]_{F_{HR}^d}^{\sim} 
  }
\end{equation*}


\subsection{The derivation functor}
\label{sec:the_derivation_functor_v2}

It is well known that discrete connection forms on $\phi$ induce
connection forms on $\phi$ (see Section 5.2
in~\cite{ar:leok_marsden_weinstein-a_discrete_theory_of_connections_on_principal_bundles}):
if $\DC\in \Sigma_C^d(\jcU)$, then $\CC:TQ\rightarrow \jgg$ defined by
\begin{equation}\label{eq:CC_from_DC-def}
  \CC(v_q)\CE  D_2\DC(q,q)(v_q)
\end{equation}
---where, as usual, $D_2$ is the differential restricted to the second
component of the Cartesian product--- is in $\Sigma_C$.

\begin{prop}
  The map $F_C:\Sigma_C^d(\jcU)\rightarrow \Sigma_C$ with
  $F_C(\DC)\CE \CC$ for $\CC$ given by~\eqref{eq:CC_from_DC-def} is
  well defined.
\end{prop}

\begin{example}\label{ex:trivial_budle:derivation}
  Let $Q\CE M\times G$ and $p_1:Q\rightarrow M$ be the trivial bundle of
  Example~\ref{ex:trivial_bundle:cc}.  Let $\DC$ be a discrete
  connection form on $p_1$ and $\CC\CE F_C(\DC)$. If
  $N:M\times M \rightarrow G$ and $\eta\in\Omega^1(M,\jgg)$ are their
  associated local expressions from Examples~\ref{ex:trivial_budle:dc}
  and~\ref{ex:trivial_bundle:cc} then
  \begin{equation*}
        \eta(\delta m) = D_2 N(m,m)(\delta m) 
  \end{equation*}
  whenever $m\in M$ and $\delta m \in T_mM$. 
\end{example}

Inspired by this idea and the Lie functor used in the Lie groupoid --
Lie algebroid context
(see~\cite{bo:mackenzie-general_theory_of_lie_groupoids_and_algebroids},
p. 125) we introduce next a functor
$\VFunctor:\FBS_M\rightarrow \VBC_M$ that will be used to assign to
each ``discrete object'' a corresponding ``continuous object''.

Let $(E,\psi,M,F)$ be a fiber bundle and $\sigma\in \Gamma(E)$, the
set of all smooth global sections of $E$. As $\psi$ is a submersion,
$\ker(T\psi)\subset TE$ is a vector bundle over $E$. Then,
$\VFunctor(E)\CE \sigma^*\ker(T\psi)$ is a vector bundle over $M$ with
rank $\dim(F)$. In fact, $p_1^r:\VFunctor(E)\rightarrow M$ ---the map
induced by $p_1:M\times TE\rightarrow M$--- is isomorphic to
$\psi\circ \tau_E|_{\ker(T\psi)}:\ker(T\psi)|_{\sigma(M)}\rightarrow
M$ as vector bundles, via the morphism
$(m,\delta e) \mapsto \delta e$. Now, for any
$f\in\hom_{\FBS_M}((E,\sigma),(E',\sigma'))$, we have that
$Tf(\ker(T\psi)) \subset \ker(T\psi')$, so that
$Tf(\ker(T\psi)|_{\sigma(M)}) \subset
\ker(T\psi')|_{\sigma'(M)}$. Thus, $Tf$ induces an element of
$\hom_{\VBC_M}(\VFunctor(E),\VFunctor(E'))$ that we denote by
$\VFunctor(f)$. It is easy to check that $\VFunctor$ is a well defined
functor from $\FBS_M$ into $\VBC_M$.

Notice that the definition of $\VFunctor(f)$ only requires that $f$ be
defined in an open neighborhood of $\sigma(M)$ in $E$, so that we can
take $\VFunctor(f)$ for any semi-local morphism~$f$.

\begin{example}\label{ex:Vfunctor_XxY}
  Let $f:X\rightarrow Y$ be a smooth map. Consider the fiber bundle
  $\psi\CE p_1:X\times Y\rightarrow X$ with the section
  $\sigma_f(x)\CE (x,f(x))$. Then, as
  $T_{(x,y)}\psi(\delta x, \delta y)=\delta x$, we have
  $\ker(T\psi)_{(x,y)} = \{0\}\oplus T_yY$, so that
  \begin{equation*}
    \VFunctor((X\times Y,\sigma_f)) = \sigma_f^*(\ker(T\psi)) =
    \{(x,(0,\delta y)):\delta y\in T_{f(x)}Y \text{ and } x\in X\}.
  \end{equation*}
  As $(x,(0,\delta y))\mapsto (x,\delta y)$ defines an isomorphism
  in $\VBC_X$, we have that
  \begin{equation*}
    \VFunctor((X\times Y,\sigma_f)) \simeq f^*TY \stext{ in }\VBC_X.
  \end{equation*}
  If, in particular, $\phi:Q\rightarrow M$ is a principal $\SG$-bundle
  and choose $f\CE \phi$, then
  $\sigma_{Q\times M}(q)\CE  \sigma_f(q) = (q,\phi(q))$ and
  \begin{equation}\label{eq:Vfunctor_QxM_iso}
    \VFunctor((Q\times M),\sigma_{(Q\times M)}) \simeq \phi^*TM \stext{ in }
    \VBC_Q.
  \end{equation}
\end{example}

\begin{example}\label{ex:Vfunctor_XxY-cases}
  Here we consider two special cases of
  Example~\ref{ex:Vfunctor_XxY}. If $X$ is any manifold and $f=id_X$,
  then $\sigma_{(X\times X)}(x)\CE \sigma_{f}(x)=(x,x)$ and
  \begin{equation}\label{eq:Vfunctor_XxX_iso}
    \VFunctor((X\times X),\sigma_{(X\times X)}) \simeq TX \stext{ in } \VBC_X.
  \end{equation}
  If $Y\CE \SG$ is a Lie group and $f:X\rightarrow \SG$ is the constant
  map $f(x)=e$, the identity element of $\SG$, then
  $\sigma_{X\times\SG}(x)\CE \sigma_f(x)=(x,e)$ and
  \begin{equation}\label{eq:Vfunctor_XxG_iso}
    \VFunctor((X\times \SG),\sigma_{(X\times \SG)}) \simeq X\times\jgg \stext{ in } \VBC_X,
  \end{equation}
  where $\jgg\CE \lie{\SG}$. 
\end{example}

\begin{example}\label{ex:Vfunctor_tiDC}
  Let $\phi:Q\rightarrow M$ be a principal $\SG$-bundle and
  $\jcU\subset Q \times Q$ be of $D$-type. Given $\DC\in
  \Sigma_C^d(\jcU)$ define
  $\ti{\DC}:Q\times Q\rightarrow Q\times \SG$ by
  $\ti{\DC}(q,q')\CE (q,\DC(q,q'))$. As
  $\ti{\DC}\in\hom_{\FBS_Q}(Q\times Q,Q\times\SG)$, using the
  characterizations provided by Example~\ref{ex:Vfunctor_XxY-cases},
  we have $\VFunctor(\ti{\DC}):TQ\rightarrow \underline{\jgg}$ (where
  $\underline{\jgg}$ is the trivial bundle over $Q$ with fiber $\jgg$)
  and, for $\delta q\in T_qQ$,
  \begin{equation*}
    T_{(q,q)}\ti{\DC}(0,\delta q) =
    (T_{(q,q)}p_1(0,\delta q),T_{(q,q)}\DC(0,\delta q)) =
    (0,D_2\DC(q,q)(\delta q)),
  \end{equation*}
  and we conclude that
  \begin{equation}\label{eq:Vfunctor_tiDC}
    \VFunctor(\ti{\DC})(\delta q) = D_2\DC(q,q)(\delta q).
  \end{equation}
\end{example}

Comparing~\eqref{eq:Vfunctor_tiDC} with~\eqref{eq:CC_from_DC-def} we
see that the assignment $\DC\mapsto \CC$ given
by~\eqref{eq:CC_from_DC-def} can be seen as coming from the functor
$\VFunctor$. In the same spirit, we look at discrete horizontal lifts.

\begin{example}\label{ex:Vfunctor_HLd}
  Let $\phi:Q\rightarrow M$ be a principal $\SG$-bundle and
  $\jcU\subset Q \times Q$ be of $D$-type. Given $\HLd{}\in
  \Sigma_H^d(\jcU)$, it follows immediately that
  $\HLd{}\in\hom_{\FBS_Q}((Q\times M,\sigma_{Q\times M}),(Q\times
  Q,\sigma_{Q\times Q}))$. As, for $(q,\delta m)\in \phi^*(TM)$,
  \begin{equation*}
    T_{(q,\phi(q))}\HLd{}(0,\delta m) = D_2\HLd{}(q,\phi(q))(\delta m),
  \end{equation*}
  the characterizations given in Examples~\ref{ex:Vfunctor_XxY}
  and~\ref{ex:Vfunctor_XxY-cases} lead to
  \begin{equation}\label{eq:Vfunctor_HLd}
    \VFunctor(\HLd{})(q,\delta m) = 
    D_2\HLds{}(q,\phi(q))(\delta m).
  \end{equation}
\end{example}

In fact, there is the following result relating discrete and
continuous horizontal lifts.
\begin{prop}\label{prop:F_H-def}
  The map $F_H:\Sigma_H^d(\jcU)\rightarrow \Sigma_H$ with
  $F_H(\HLd{}) \CE  \VFunctor(\HLd{})$ (computed
  in~\eqref{eq:Vfunctor_HLd}) is well defined.
\end{prop}

Last, we apply the $\VFunctor$ functor to right splittings
of~\eqref{eq:DAS-def} in an attempt to find right splittings
of~\eqref{eq:AS-def}.

\begin{example}\label{ex:Vfunctor_(QxQ)/G}
  Let $\phi:Q\rightarrow M$ be a principal $\SG$-bundle. Consider the
  fiber bundle
  $\psi\CE \widecheck{\phi\circ p_1}:(Q\times Q)/\SG\rightarrow M$ with
  the section $\sigma_{(Q\times Q)/\SG}$ defined
  in~\eqref{eq:DAS_spaces-sigma-def}. Then, as
  \begin{equation*}
    \begin{split}
      T_{[(q,q')]}\psi(T_{(q,q')}\pi^{Q\times Q,\SG}(\delta q,\delta
      q')) =& T_{(q,q')}(\psi\circ \pi^{Q\times Q,\SG})(\delta
      q,\delta q')) \\=& T_{(q,q')}(\phi\circ p_1)(\delta q,\delta q')) =
      T_q\phi(\delta q),
    \end{split}
  \end{equation*}
  we have that
  \begin{equation*}
    \begin{split}
      \ker(T\psi)_{[(q,q')]} =& \{T_{(q_0,q_0')}\pi^{Q\times
        Q,\SG}(\delta q_0, \delta q_0'):\delta q_0\in {\mathcal
        V}^Q_{q_0}, \delta q_0'\in T_{q_0'}Q, [(q,q')]=[(q_0,q_0')]\}
      \\=& \{T_{(q_0,q_0')}\pi^{Q\times Q,\SG}(0, \delta q_0''):\delta
      q_0''\in T_{q_0'}Q, [(q,q')]=[(q_0,q_0')]\}.
    \end{split}
  \end{equation*}
  Consequently,
  \begin{equation*}
    \begin{split}
      \VFunctor((Q\times Q)/\SG) =& \{(m,T_{(q,q)}\pi^{Q\times
        Q,\SG}(0,\delta q'')):\phi(q)=m, \delta q''\in T_qQ\}.
    \end{split}
  \end{equation*}
  In addition, mapping
  $(m,T_{(q,q)}\pi^{Q\times Q,\SG}(0,\delta q'')) \mapsto
  \pi^{TQ,\SG}(\delta q'')$, defines an isomorphism in $\VBC_M$, so
  that
  \begin{equation}\label{eq:Vfunctor_(QxQ)/G_iso}
    \VFunctor((Q\times Q)/\SG) \simeq (TQ)/\SG \stext{ in } \VBC_M.
  \end{equation}
\end{example}

\begin{example}\label{ex:Vfunctor_F2}
  Let $F_2$ be the map defined in~\eqref{eq:DAS-def}
  and~\eqref{eq:DAS_maps-def}. As $F_2$ is a morphism in the $\FBS_M$
  category, we may compute $\VFunctor(F_2)$. Recalling the
  characterizations of $\VFunctor(M\times M)$ and
  $\VFunctor((Q\times Q)/\SG)$ given in the
  Examples~\ref{ex:Vfunctor_XxY-cases} and~\ref{ex:Vfunctor_(QxQ)/G}, we
  have that, for any $\delta q''\in T_qQ$,
  \begin{equation*}
    \begin{split}
      T_{[(q,q)]}F_2(T_{(q,q)}\pi^{Q\times Q,\SG}(0,\delta q'')) =&
      T_{(q,q)}(F_2\circ \pi^{Q\times Q,\SG})(0,\delta q'') \\=&
      T_{(q,q)}(\phi\times \phi)(0,\delta q'') = (0,T_q\phi(\delta
      q'')).
    \end{split}
  \end{equation*}
  Taking into account~\eqref{eq:Vfunctor_XxX_iso}
  and~\eqref{eq:Vfunctor_(QxQ)/G_iso}, this leads to
  \begin{equation*}
    \VFunctor(F_2)([\delta q'']) = T_q\phi(\delta q'').
  \end{equation*}
  In other words,
  \begin{equation*}
    \VFunctor(F_2) = \phi_2
  \end{equation*}
  for the map $\phi_2$ appearing in~\eqref{eq:AS-def}.
\end{example}

\begin{remark}
  It can also be seen that $\VFunctor(\ti{\SG}) = \ti{\jgg}$ and
  $\VFunctor(F_1) = \phi_1$, so that the $\VFunctor$ functor maps the
  Discrete Atiyah Sequence~\eqref{eq:DAS_maps-def} into the Atiyah
  Sequence~\eqref{eq:AS-def}.
\end{remark}

\begin{prop}
  The map $F_R:\Sigma_R^d(\jcU)\rightarrow \Sigma_R$ defined by
  $F_R(s_d)\CE \VFunctor(s_d)$ is well defined.
\end{prop}

\begin{proof}
  We need to check that, for $s_d\in \Sigma_R^d(\jcU)$,
  $s\CE \VFunctor(s_d)$ is a right splitting
  of~\eqref{eq:AS-def}. Using the characterizations given by
  Examples~\ref{ex:Vfunctor_XxY-cases}, \ref{ex:Vfunctor_(QxQ)/G}
  and~\ref{ex:Vfunctor_F2}, we have that
  $s\in\hom_{\VBC_M}(TM,(TQ)/\SG)$ and
  \begin{equation*}
    \phi_2\circ s = \VFunctor(F_2)\circ \VFunctor(s_d) = \VFunctor(F_2\circ s_d) = \VFunctor(id_{\jcU''}) = id_{TM},
  \end{equation*}
  concluding the argument.
\end{proof}

We can summarize the different sets and maps that we have constructed
so far via the following diagram in the category of sets,
\begin{equation}\label{eq:full_D_and_C_diagram}
  \begin{split}
    \xymatrix{ {\Sigma_C^d(\jcU)} \ar[rr]^{F_{CH}^d}
      \ar[dr]_{F_{CR}^d}
      \ar@{-->}[ddd]_{F_C}& {} & {\Sigma_H^d(\jcU)}
      \ar[dl]^{F_{HR}^d} \ar@{-->}[ddd]^{F_H}\\
      {} & {\Sigma_R^d(\jcU)} \ar@{-->}[d]^{F_R} & {} \\
      {} & {\Sigma_R} & {} \\
      {\Sigma_C} \ar[rr]_{F_{CH}} \ar[ur]^{F_{CR}} & {} &
      {\Sigma_H} \ar[lu]_{F_{HR}} }
  \end{split}
\end{equation} 
where the full arrows are bijections. Careful unraveling of the
definitions leads to the following result.

\begin{prop}\label{prop:commutativity_of_full_D_and_C_diagram}
  Diagram~\eqref{eq:full_D_and_C_diagram} in the \setC category is
  commutative.
\end{prop}


\section{The flat case}
\label{sec:the_flat_case}

In this section we review how flat connections and flat discrete
connections may be viewed as splittings of sequences in the categories
of Lie algebroids and local Lie groupoids. Then, using an
integration result for morphisms in the category of Lie algebroids,
we prove that the Integration Problem in the flat case admits a unique
solution.


\subsection{Lie groupoids and Lie algebroids}
\label{sec:lie_groupoids_and_lie_algebroids}

We start by reviewing the basic notions of local Lie groupoid and Lie
algebroid, mostly to fix the notation. There are several slightly
different notions of \emph{local} Lie groupoid. We essentially follow
the one introduced in p.40
of~\cite{ar:coste_dazord_weinstein-grupoides_symplectiques}, that is
also used
in~\cite{ar:fernandez_juchani_zuccalli-discrete_connections_on_principal_bundles_the_discrete_atiyah_sequence}.

\begin{definition}\label{def:local_lie_groupoid}
  A \jdef{local Lie groupoid} consists of smooth manifolds $G$ and $M$
  together with submersions $\alpha,\beta:G\rightarrow M$ as well as a
  diffeomorphism $i:G\rightarrow G$, and smooth maps
  $\epsilon:M\rightarrow G$ and $m:G_m\rightarrow G$, where
  $G_m\subset G_2\CE G\FP{\beta}{\alpha} G$ is an open subset, all
  subject to the conditions stated below. It is convenient to write
  $g_1g_2\CE m(g_1,g_2)$ and $g^{-1}\CE i(g)$.
  \begin{enumerate}
  \item \label{it:lLgpd-epsilon}
    $\alpha\circ \epsilon = id_M = \beta\circ \epsilon$.
  \item \label{it:lLgpd-inverse} For all $(g_1,g_2)\in G_m$, we have
    that $(g_2^{-1}, g_1^{-1}) \in G_m$ and
    $g_2^{-1}g_1^{-1} = (g_1g_2)^{-1}$.
  \item \label{it:lLgpd-identity} For all $g\in G$ we have that
    $(\epsilon(\alpha(g)),g), (g,\epsilon(\beta(g))) \in G_m$ and
    $\epsilon(\alpha(g)) g = g = g \epsilon(\beta(g))$.
  \item \label{it:lLgpd-inverses} For all $g\in G$ we have that
    $(g,g^{-1}),(g^{-1},g)\in G_m$ and $gg^{-1} = \epsilon(\alpha(g))$
    and $g^{-1}g=\epsilon(\beta(g))$.
  \item \label{it:lLgpd-associativity} If
    $(g_1,g_2), (g_2,g_3), (g_1,g_2g_3)\in G_m$, then
    $(g_1g_2,g_3) \in G_m$ and $(g_1g_2)g_3 = g_1(g_2g_3)$.
  \end{enumerate}
  We denote a local Lie groupoid as above by $G\gpoidarrows M$. If
  $G_m=G_2$ the local Lie groupoid is a \jdef{Lie groupoid}\footnote{Our
    definition of Lie groupoid is, in a sense, opposed to Definition
    1.1.1
    of~\cite{bo:mackenzie-general_theory_of_lie_groupoids_and_algebroids}. For
    example, in Definition~\ref{def:local_lie_groupoid}, the
    ``multiplicable elements'' are $G_2 = G\FP{\beta}{\alpha} G$,
    where
    in~\cite{bo:mackenzie-general_theory_of_lie_groupoids_and_algebroids}
    the analogous set is $G_2 = G\FP{\alpha}{\beta} G$. Of course, the
    results are the same, for both conventions, albeit with some
    cosmetic rewriting.}.
\end{definition}

\begin{definition}\label{def:morphism_of_local_lie_groupoid}
  Let $G\gpoidarrows M$ and $G'\gpoidarrows M'$ be local Lie
  groupoids. A smooth map $F:G\rightarrow G'$ is a \jdef{morphism of
    local Lie groupoids} if
  \begin{enumerate}
  \item $F(\epsilon_G(M))\subset \epsilon_{G'}(M')$, and
  \item for all $(g_1,g_2)\in G_m$, we have that
    $(F(g_1),F(g_2))\in G'_m$ and $F(g_1g_2)=F(g_1)F(g_2)$.
  \end{enumerate}
\end{definition}

The local Lie groupoids together with their
morphisms form a category that we denote by $\lLgpdC$. In what
follows, we will be mostly interested in the subcategory $\lLgpdC_M$
of local Lie groupoids over a smooth manifold $M$ and whose morphisms cover
$id_M$.

\begin{example}\label{ex:DAS_is_sequence_in_lLgpdC}
  The Discrete Atiyah Sequence~\eqref{eq:DAS-def} is a sequence in the
  $\lLgpdC_M$ category. In fact, it is a sequence of Lie
  groupoids. Indeed, $M\times M$ is what is known as the \jdef{pair
    groupoid} and the structure of $(Q\times Q)/\SG$ is induced by
  that of the pair groupoid $Q\times Q$. The structure maps of
  $\ti{\SG}$ are $\alpha=\beta=$ the fiber bundle projection on $M$,
  $m((q,g_0),(q,g_1))\CE (q,g_0g_1)$, $\epsilon(m)\CE \sigma_{\ti{\SG}}(m)$
  for $\sigma_{\ti{\SG}}$ defined in~\eqref{eq:DAS_spaces-sigma-def}
  and $i(q,g)\CE (q,g^{-1})$. For more details, see Section 5.1
  in~\cite{ar:fernandez_juchani_zuccalli-discrete_connections_on_principal_bundles_the_discrete_atiyah_sequence}.
\end{example}

\begin{example}
  Let $\phi:Q\rightarrow Q/\SG$ be a principal $\SG$-bundle and
  $\jcU\subset Q\times Q$ be of symmetric $D$-type. Then we have the
  restriction of the Discrete Atiyah Sequence~\eqref{eq:DAS-def} to
  $\jcU$,
  \begin{equation}\label{eq:DAS_in_lLgpd}
    \xymatrix{
      {\ti{\SG}} \ar[r]^-{F_1} & {\jcU/\SG} \ar[r]^-{F_2} & {\jcU''}
    }
  \end{equation}
  where the maps $F_1$ and $F_2$ are the appropriate restriction and
  co-restrictions of the maps defined
  in~\eqref{eq:DAS_maps-def}. Then,~\eqref{eq:DAS_in_lLgpd} is a
  sequence in the $\lLgpdC_M$ category. For more details, see Section
  5.2
  in~\cite{ar:fernandez_juchani_zuccalli-discrete_connections_on_principal_bundles_the_discrete_atiyah_sequence}.
\end{example}

\begin{definition}\label{def:lie_algebroid}
  A \jdef{Lie algebroid} on the smooth manifold $M$ is triple
  $(A,[,],\rho)$ where $A\in\VBC_M$, $[,]$ is a Lie bracket on
  $\Gamma(A)$ and $\rho\in\hom_{\VBC_M}(A,TM)$ is the \jdef{anchor
    map}, satisfying
  \begin{equation*}
    [X,fY]=(\rho(X)f)\, Y+ f[X,Y], \stext{ for all }
    X,Y\in \Gamma(A) \text{ and }
    f\in C^\infty(M).
  \end{equation*}
\end{definition}

\begin{definition}\label{def:morphism_of_lie_algebroids}
  If $(A_1,[,]_1,\rho_1)$ and $(A_2,[,]_2,\rho_2)$ are two Lie
  algebroids over $M$, a \jdef{morphism of Lie algebroids} is a map
  $F\in \hom_{\VBC_M}(A_1,A_2)$ such that
  \begin{enumerate}
  \item $F([X,Y]_1) = [F(X),F(Y)]_2$ for every $X,Y\in \Gamma(A_1)$, and
  \item $\rho_2\circ F = \rho_1$.
  \end{enumerate}
\end{definition}

The Lie algebroids over a smooth manifold $M$ together with their
morphisms form a category that we denote by $\LalgbdC_M$.

\begin{example}\label{ex:AS_is_sequence_in_LalgbdpdC}
  The Atiyah Sequence~\eqref{eq:AS-def} is a sequence in the
  $\LalgbdC_M$ category. For the details, see Section 3.2
  in~\cite{bo:mackenzie-general_theory_of_lie_groupoids_and_algebroids}.
\end{example}

The classical \jdef{Lie functor} assigns all Lie groups and morphisms,
the corresponding Lie algebras and morphisms, via a derivation
process. This functor has been extended to the categories of (local)
Lie groupoids and Lie algebroids
(see~\cite{bo:mackenzie-general_theory_of_lie_groupoids_and_algebroids},
p. 125). Because of this construction it is customary to say that the
morphism of Lie algebroids obtained via the Lie functor from a
morphism of Lie groupoids is ``its derivative'' and, conversely, call
the latter ``an integral'' of the former. Notice that, for
(transitive) Lie groupoids, forgetting the brackets and anchor map,
the Lie functor becomes the $\VFunctor$ functor introduced in
Section~\ref{sec:the_derivation_functor_v2}. Also, it can be checked
that applying the Lie functor from $\lLgpdC_M$ into $\LalgbdC_M$ to
the Discrete Atiyah Sequence~\eqref{eq:DAS-def} of a principal bundle
gives its Atiyah Sequence~\eqref{eq:AS-def}.

The following existence and uniqueness result of A. Cabrera et al. is
part of a version of Lie's Second Theorem in the context of Lie
groupoids. It is the key to solving the Integration Problem for flat
connections.

\begin{theorem}\label{th:integration_of_morphisms_of_lie_algebroids}
  Let $G_1\gpoidarrows M$ and $G_2\gpoidarrows M$ be local Lie
  groupoids with Lie algebroids $A_1$ and $A_2$ respectively, and let
  $f:A_1\rightarrow A_2$ be a Lie algebroid morphism. There exists a
  local Lie groupoid morphism $F:G_1\rightarrow G_2$ integrating $f$,
  and any two such integrations coincide around the unit section.
\end{theorem}

\begin{proof}
  Theorem 2.4 on p. 7
  of~\cite{ar:cabrera_marcut_salazar-on_local_integration_of_lie_brackets}.
\end{proof}

\begin{remark}
  The notion of local Lie groupoid used
  in~\cite{ar:cabrera_marcut_salazar-on_local_integration_of_lie_brackets}
  is slightly different from the one that we are using in this
  paper. Still, as the authors prove in Proposition 2.1
  of~\cite{ar:cabrera_marcut_salazar-on_local_integration_of_lie_brackets}
  the germ of any local Lie groupoid in their sense can be represented
  by a local Lie groupoid in our sense. Also, the meaning of
  $F:G_1\rightarrow G_2$ being a ``local Lie groupoid morphism''
  in~\cite{ar:cabrera_marcut_salazar-on_local_integration_of_lie_brackets}
  is that there exist open subsets $U\subset G_1$ and
  $W\subset (G_1)_2$ such that $\epsilon_1(M)\subset U$ and
  $\epsilon_1(M)_2\CE (\epsilon_1(M))^2\cap (G_1)_2\subset W$, that
  $F:U\rightarrow G_2$ is smooth, that
  $F(\epsilon_1(M)) \subset \epsilon_2(M)$ and that for all
  $(g,g')\in W$, $F(gg') = F(g)F(g')$.
\end{remark}

\begin{corollary}\label{cor:integration_of_morphs_lLalgbd_and_lLgpd}
  Let $G_1\gpoidarrows M$ and $G_2\gpoidarrows M$ be Lie groupoids
  over $M$ with associated Lie algebroids $A_1\rightarrow M$ and
  $A_2\rightarrow M$. If $f:A_1\rightarrow A_2$ is a morphism of Lie
  algebroids, then, there are open subsets $U\subset G_1$ and
  $U_m\subset (G_1)_2$ such that
  \begin{enumerate}
  \item $U$ with the structure maps induced by $G_1$ is a local Lie
    groupoid whose multiplication $m_U$ is defined on $U_m$.
  \item There exists $F\in\hom_{\lLgpdC_M}(U,G_2)$ such that
    $\LieFunctor(F)=f$.
  \end{enumerate}
\end{corollary}

\begin{proof}
  By Theorem~\ref{th:integration_of_morphisms_of_lie_algebroids}, there are
  open subsets $V\subset G_1$ and $W\subset (G_1)_2$ such that
  $\epsilon_1(M)\subset V$ and $\epsilon_1(M)_2\subset W$ and there exists a
  smooth map $F:V\rightarrow G_2$ such that $F(\epsilon_1(M))\subset
  \epsilon_2(M)$, and, for $(g,g')\in W$, $F(gg') =F(g)F(g')$ (multiplications
  in $G_1$ and $G_2$). Furthermore, $\LieFunctor(F)=f$.

  There is a claim in the proof of Proposition 2.1
  in~\cite{ar:cabrera_marcut_salazar-on_local_integration_of_lie_brackets}
  (but see Remark~\ref{rem:seeapp} later on) to the effect that the
  open subsets $U^{(2)}\CE \left(U\times U\right)\cap (G_1)_2$ for
  $U\subset G_1$ such that $\epsilon_1(M)\subset U$ form a basis of
  open neighborhoods of $\epsilon_1(M)_{2}\subset (G_1)_2$. Using this
  claim applied to $W\subset (G_1)_2$, we have that there is an open
  subset $V'\subset G_1$ containing $\epsilon_1(M)$ and such that
  $\epsilon_1(M)_2\subset (V')^{(2)}\subset W$.
  
  Now let $U'\CE V\cap V'$ and $U\CE U'\cap (U')^{-1}$. Then,
  $U^{-1}\subset U$, so that it inherits a local Lie groupoid over $M$
  structure from $G_1$ using Lemma 5.11
  in~\cite{ar:fernandez_juchani_zuccalli-discrete_connections_on_principal_bundles_the_discrete_atiyah_sequence};
  notice that in this construction, the multiplication map $m_U$ is
  defined on $U_m\CE U_2 \cap m_1^{-1}(U) \subset U^{(2)} \subset W$,
  so that the map $F$ that integrates $f$ is, indeed, in
  $\hom_{\lLgpdC_M}(U,G_2)$.
\end{proof}


\subsection{Flat connections in the Lie context}
\label{sec:flat_connections_in_the_lie_context}

Let $\CC$ be a principal connection on the principal $\SG$-bundle
$\phi:Q\rightarrow M$. If we let $s\CE F_{CR}(\CC)$ be the
representation of the connection as a right splitting of the Atiyah
Sequence~\eqref{eq:AS-def}, we see that, for any $X,Y\in\VF(M)$,
$s([X,Y])-[s(X),s(Y)] \in \ker(\phi_2)$. Due to the exactness of the
Atiyah Sequence, there is a unique $b_{XY}\in \ti{\jgg}$ such that
$\phi_1(b_{XY}) = s([X,Y])-[s(X),s(Y)]$.

\begin{definition}\label{def:CB}
  In the previous setting, the morphism of vector bundles over $M$,
  $\CB:TM\oplus TM\rightarrow \ti{\jgg}$ defined on vector fields by
  \begin{equation*}
    \CB(X,Y)\CE b_{XY}
  \end{equation*}
  is the \jdef{curvature} of $\CC$. Also, $\CC$ is said to be
  \jdef{flat} if $\CB=0$.
\end{definition}

\begin{remark}
  Definition~\ref{def:CB} ---Definition 5.3.8
  in~\cite{bo:mackenzie-general_theory_of_lie_groupoids_and_algebroids}---
  is not the most common definition of this notion. Still, it is the
  most appropriate for our work. Furthermore, by Proposition 5.3.14
  in~\cite{bo:mackenzie-general_theory_of_lie_groupoids_and_algebroids},
  it is equivalent to the standard one, as it appears in, for
  instance,~\cite{bo:kobayashi_nomizu-foundations-v1}.
\end{remark}

\begin{example}\label{ex:trivial_bundle:flat:conti}
  Let $M$ be a manifold and $G$ be an \emph{abelian} Lie group. Following
  Example~\ref{ex:trivial_bundle:cc}, consider the principal
  $G$-bundle $p_1:Q\CE M\times G \rightarrow M$.  Given a connection form
  $\CC$ on $p_1$ we define $\omega\in\Omega^1(M,\jgg)$ to be its local
  expression.  It is straightforward to see that $\CC$ is flat if and
  only if the differential form $\omega$ is closed.
\end{example}

Definition~\ref{def:CB} states explicitly that the curvature of a
principal connection $\CC$ is the obstruction for its associated
section $s$ to induce a morphism of Lie algebras between the
corresponding spaces of sections. An immediate consequence is the
following result.

\begin{prop}\label{prop:CC_flat_connection_iff_s_morph_in_Lalgbd}
  Let $\CC$ be a principal connection on the principal $\SG$-bundle
  $\phi:Q\rightarrow M$ and let $s\CE F_{CR}(\CC)$ be the
  representation of the connection as a right splitting of the Atiyah
  Sequence~\eqref{eq:AS-def} in the $\VBC_M$ category. Then $s$ is a
  morphism of Lie algebroids if and only if $\CC$ is flat.
\end{prop}

\begin{proof}
  By construction, $s\in \Sigma_R$ is a morphism of smooth vector
  bundles and it preserves the corresponding anchor maps. Then, as
  $\phi_1$ is injective, $s$ is a morphism of the Lie algebras of
  sections if and only if $\CB$ vanishes.
\end{proof}

There is, also, a notion of curvature for discrete connections, as we
recall next.

\begin{definition}\label{def:discrete_curvature}
  Let $\phi:Q\rightarrow M$ be a principal $\SG$-bundle,
  $\jcU\subset Q\times Q$ be of $D$-type and
  $\DC\in\Sigma_C^d(\jcU)$. Let
  \begin{equation}\label{eq:U^{(3)}-def}
    \jcU^{(3)}
    \CE  \{(q_0,q_1,q_2) \in Q^3 : (q_i,q_j)\in \jcU 
    \text{ for all } 0\leq i<j\leq 2\}.
  \end{equation}
  We define the \jdef{discrete curvature} of $\DC$ as
  $\DB: \jcU^{(3)} \rightarrow \SG$ by
  \begin{equation}\label{eq:BD-def}
    \DB(q_0,q_1,q_2) \CE  \DC(q_0,q_2)^{-1}\DC(q_1,q_2)\DC(q_0,q_1).
  \end{equation}
  We say that $\DC$ is \jdef{flat} if $\DB=e$ on $\jcU^{(3)}$.
\end{definition}

The following result shows that there is a bijection between the flat
discrete connections and the right splittings of the Discrete Atiyah
Sequence in the $\lLgpdC_M$ category.

\begin{prop}\label{prop:bijection_right_splittings_DAS_flat_DC-lLgpd}
  Let $\phi:Q\rightarrow M$ be a principal $\SG$-bundle and
  $\jcU\subset Q\times Q$ be of symmetric $D$-type. Then, the
  bijection $F_{CR}^d:\Sigma_C^d(\jcU)\rightarrow \Sigma_R^d(\jcU)$
  determines a bijection between the subset
  $\Sigma_C^{d,e}(\jcU)\subset \Sigma_C^d(\jcU)$ of flat discrete
  connections and the subset
  $\ti{\Sigma_R^d}(\jcU)\subset \Sigma_R^d(\jcU)$ of right splittings
  of extension~\eqref{eq:DAS_in_lLgpd} in the $\lLgpdC_M$ category.
\end{prop}

\begin{proof}
  Proposition 5.30
  in~\cite{ar:fernandez_juchani_zuccalli-discrete_connections_on_principal_bundles_the_discrete_atiyah_sequence}.
\end{proof}

\begin{remark}\label{rem:MECs_and_DCs}
  Other notions inspired by the connections on a principal
  $\SG$-bundle have been considered in the literature. For example,
  the multiplicative Ehresmann connections (MEC) introduced
  in~\cite{ar:fernandes_marcut-multiplicative_ehresmann_connections},
  based on the notion of connection for Lie groupoid extensions
  of~\cite{laurentGengoux_stienon_xu-non_abelian_gerbes-2009}, is a
  notion of connection for a (surjective and submersive) morphism of
  Lie groupoids $F:\jcG\rightarrow \jcH$. Essentially, it is a (not
  necessarily integrable) distribution $E$ on $\jcG$ such that
  $T\jcG=E\oplus \ker(TF)$ and that $E$ is a Lie subgroupoid of
  $T\jcG$. The case that could be compared to the setting of discrete
  connections on $\phi:Q\rightarrow M$ is that where
  $\jcG\CE (Q\times Q)/\SG$ with structure maps
  $\alpha_\jcG\CE \check{p}_1$ and $\beta_\jcG\CE \check{p}_2$, $\jcH$
  is the pair groupoid $M\times M$ and
  $F\CE (\alpha_\jcG,\beta_\jcG)$.

  Any discrete connection $\DC$ with domain $\jcU\subset Q\times Q$ on
  $\phi$ defines the semi-local splitting
  $s_R\CE F_{CR}^d(\DC) \in \Sigma_R^d(\jcU)$ (see
  Section~\ref{sec:DC_on_principal_bundles}). In general, for
  $\jcU''\CE (\phi\times \phi)(\jcU)$, the map
  $s_R:\jcU''\rightarrow (Q\times Q)/\SG$ is an embedding so that
  $\check{H}\CE s_R(\jcU'')$ is an embedded submanifold; thus,
  $E\CE T\check{H}\subset T((Q\times Q)/\SG)|_{\check{H}}$ is a
  distribution over $\check{H}$ (that, by construction, is
  integrable). It can be seen that
  $T((Q\times Q)/\SG)|_{\check{H}} = E \oplus \ker
  T(\check{p}_1,\check{p}_2)|_{\check{H}}$. Neither $\check{H}$ nor
  $E$ are usually local Lie groupoids. Still, when $\DC$ is flat,
  $s_R$ is a morphism of local Lie groupoids, $\check{H}$ is a local
  Lie groupoid and $E$ is a local Lie subgroupoid of
  $T((Q\times Q)/\SG)|_{\check{H}}$.

  So, discrete connections share some of the properties with MECs but
  are clearly different. Additional investigation may lead to further
  clarification of their relationship.
\end{remark}

The derivation process preserves the flatness condition, as we see next.

\begin{prop}\label{prop:F_C_preserves_flatness}
  Let $\phi:Q\rightarrow M$ be a principal $\SG$-bundle and
  $\DC\in \Sigma_C^{d,e}(\jcU)$ for a $D$-type open subset
  $\jcU\subset Q\times Q$. Then, $F_C(\DC) \in \Sigma_C$ is a flat
  connection in $\phi$.
\end{prop}

\begin{proof}
  As $\DC\in \Sigma_C^{d,e}(\jcU)$, by
  Proposition~\ref{prop:bijection_right_splittings_DAS_flat_DC-lLgpd},
  $s_d\CE F_{CR}^d(\DC) \in \ti{\Sigma_R^d}(\jcU)$. Then,
  $\xymatrix{{(Q\times Q)/\SG} \ar[r]_-{F_2} & {M\times M}
    \ar@/_1pc/[l]_{s_d}}$ is a diagram in the $\lLgpdC_M$ category and
  $F_2\circ s_d =id_{M\times M}$. Applying the Lie functor to the
  previous diagram and recalling Example~\ref{ex:Vfunctor_F2} we
  obtain the diagram
  $\xymatrix{{(TQ)/\SG} \ar[r]_-{\phi_2} & {TM} \ar@/_1pc/[l]_{s}}$ in
  the $\LalgbdC_M$ category, where
  $s\CE \LieFunctor(s_d) = \VFunctor(s_d)$. Applying the same procedure
  to the identity, we obtain $\phi_2 \circ s = id_{TM}$, proving that
  $s$ is a right splitting of~\eqref{eq:AS-def} in the $\LalgbdC_M$
  category. Therefore by
  Proposition~\ref{prop:CC_flat_connection_iff_s_morph_in_Lalgbd},
  $\CC\CE F_{RC}(s)$ is flat. Using
  Proposition~\ref{prop:commutativity_of_full_D_and_C_diagram} the
  statement follows because
  $\CC = F_{RC}(s) =F_{RC}(\VFunctor(s_d)) = F_{RC}(F_R(s_d)) =
  F_{RC}(F_R(F_{CR}^d(\DC))) = F_C(\DC)$.
\end{proof}

The next result puts together the interpretation of flat connections
as morphisms of Lie groupoids and algebroids, turning the Integration
Problem into one of integration of morphisms of Lie algebroids, that
can be solved using
Theorem~\ref{th:integration_of_morphisms_of_lie_algebroids}.

\begin{theorem}\label{th:solution_of_IP-flat_case}
  Let $\phi:Q\rightarrow M$ be a principal $\SG$-bundle and
  $\CC\in \Sigma_C$ that is flat. Then,
  \begin{enumerate}
  \item \label{it:solution_of_IP-flat_case-existence} there is a
    subset $\jcU\subset Q\times Q$ of symmetric $D$-type and
    $\DC\in \Sigma_C^d(\jcU)$ flat that integrates $\CC$, that is,
    $\CC=F_C(\DC)$.
  \item \label{it:solution_of_IP-flat_case-uniqueness} If
    $\jcU_1,\jcU_2$ are symmetric of $D$-type such that there are
    $\DCp{1} \in \Sigma_C^d(\jcU_1)$ and
    $\DCp{2} \in \Sigma_C^d(\jcU_2)$, both of which are flat and
    integrate $\CC$, then there is a subset
    $\jcU\subset \jcU_1\cap\jcU_2$ of symmetric $D$-type set such that
    $\DCp{1}$ and $\DCp{2}$ agree on $\jcU$.
  \end{enumerate}
\end{theorem}

\begin{proof}
  Recall that $M\times M$ and $(Q\times Q)/\SG$ are Lie groupoids over
  $M$ whose corresponding Lie algebroids are $TM$ and $(TQ)/\SG$ (see
  Example 3.5.11 and Theorem 4.5.7
  in~\cite{bo:mackenzie-general_theory_of_lie_groupoids_and_algebroids}).

  Let $s\CE F_{CR}(\CC) \in \Sigma_R$; as $\CC$ is flat, by
  Proposition~\ref{prop:CC_flat_connection_iff_s_morph_in_Lalgbd},
  $s:TM\rightarrow (TQ)/\SG$ is a morphism of Lie algebroids over
  $M$. Then, by
  Corollary~\ref{cor:integration_of_morphs_lLalgbd_and_lLgpd}, there
  are open subsets $U\subset M\times M$ and $U_m\subset (M\times M)_2$
  such that, with the structure induced by that of the Lie groupoid
  $M\times M$ and multiplication defined on $U_m$, $U$ is a local Lie
  groupoid over $M$ and there is
  $s_d\in\hom_{\lLgpdC_M}(U,(Q\times Q)/\SG)$ such that
  $\LieFunctor(s_d)=s$, that is, $s_d$ integrates $s$. Notice that
  $H\CE F_2\circ s_d \in \hom_{\lLgpdC_M}(U,M\times M)$ satisfies
  $\lie{H} = \lie{F_2 \circ s_d} = \lie{F_2}\circ \lie{s_d} =
  \phi_2\circ s = id_{TM} = \lie{id_U}$. Then, by the uniqueness part
  of Theorem~\ref{th:integration_of_morphisms_of_lie_algebroids},
  shrinking $U$ if necessary, we have that $H=id_U$, so that $s_d$ is
  a section of $F_2$.
  Let $\jcU\CE (\phi\times \phi)^{-1}(U)\subset Q\times Q$. It is easy
  to check that $\jcU$ is of symmetric $D$-type and $\jcU''=U$. Hence,
  $s_d\in\hom_{\lLgpdC_M}(U,(Q\times Q)/\SG)$ implies that
  $\DC\CE (F_{CR}^{d})^{-1}(s_d) \in \Sigma_C^{d,e}(\jcU)$. Following
  the arrows in diagram~\eqref{eq:full_D_and_C_diagram}, we have that
  $F_C(\DC) = \CC$, proving
  point~\eqref{it:solution_of_IP-flat_case-existence} of the
  statement.

  Given $\DCp{j}\in \Sigma_C^d(\jcU_j)$, $j=1,2$, both flat and
  integrating $\CC$, we define
  $s_d^j\CE F_{CR}^d(\DCp{j}) \in \ti{\Sigma_R^d}(\jcU_j)$ and
  $s\CE F_{CR}(\CC)\in \Sigma_R$; as $\CC$ is flat, by
  Proposition~\ref{prop:CC_flat_connection_iff_s_morph_in_Lalgbd},
  $s:TM\rightarrow (TQ)/\SG$ is a morphism of Lie algebroids over $M$
  and $s = \LieFunctor(s_d^j)$ for $j=1,2$. Then, by
  Theorem~\ref{th:integration_of_morphisms_of_lie_algebroids}, there
  is an open subset $V\subset \jcU_1''\cap \jcU_2''\subset M\times M$
  containing $\Delta_M$ such that $s_d^1|_{V}=s_d^2|_{V}$. Defining
  $V'\CE V\cap V^{-1}$ (the inverse in the pair groupoid structure of
  $M\times M$), we see that $\Delta_M\subset V'\subset V$, so that
  $\jcU\CE (\phi\times \phi)^{-1}(V')$ is symmetric of $D$-type and
  $\jcU''=V'$. Last, $s_d^1|_{V} = s_d^2|_{V}$ implies that
  $\DCp{1}|_{\jcU} = \DCp{2}|_{\jcU}$, concluding the proof of
  point~\eqref{it:solution_of_IP-flat_case-uniqueness} in the
  statement.
\end{proof}

\begin{remark}\label{rem:seeapp}
  Corollary~\ref{cor:integration_of_morphs_lLalgbd_and_lLgpd},
  referenced in the proof of Theorem~\ref{th:solution_of_IP-flat_case}
  above, relies on a claim made in the proof of Proposition~2.1
  of~\cite{ar:cabrera_marcut_salazar-on_local_integration_of_lie_brackets}. Unfortunately,
  the proof of that claim is only slightly hinted at in the paper.  We
  provide an alternative, more detailed, explanation of the part of
  the claim that is needed in the proof of
  Corollary~\ref{cor:integration_of_morphs_lLalgbd_and_lLgpd} in
  Proposition~\ref{prop:basis}.
\end{remark}


\section{Existence of solution in the general case}
\label{sec:existence_of_solution_in_the_general_case}

In this section we abandon the special case of flat connections and
return to the arbitrary curvature case. We will be able to integrate
arbitrary principal connections provided that we have the additional
data of a smooth retraction. Then we state sufficient conditions on
the structure group of the principal bundle under which these
retractions do exist.


\subsection{Smooth retractions}
\label{sec:smooth_retractions}

Following~\cite[Definition
4.1.1]{bo:absil_mahony_sepulchre-optimization_algorithms_on_matrix_manifolds}
we introduce the following notion.

\begin{definition}\label{def:retraction}
  Let $X$ be a smooth manifold and denote by $Z_X$ the image of its
  zero section in the tangent bundle.  A \jdef{smooth retraction} on
  $X$ is a smooth map $R:\jcZ\rightarrow X$, where $\jcZ\subset TX$ is
  an open neighborhood of $Z_X$, such that, if we denote
  $R_x\CE R|_{T_xX\cap \jcZ}$ for any $x\in X$, the following properties
  are valid.
  \begin{enumerate}
  \item \label{it:retraction_0} $R_x(0_x)=x$ for all $x\in X$, where
    $0_x$ is the null vector in $T_xX$.
  \item \label{it:retraction_rigidity} Identifying canonically
    $T_{0_x}T_xX \simeq T_xX$, we have $T_{0_x}R_x = id_{T_xX}$.
  \end{enumerate}
\end{definition}

\begin{remark}
  Condition~\eqref{it:retraction_rigidity} in
  Definition~\ref{def:retraction} is equivalent to requiring that
  \begin{equation}\label{eq:retraction_cond_2_explicit}
    \frac{d}{dt}\bigg|_{t=0} R_x(tv_x) = v_x \in T_xX
    \stext{ for all } v_x\in T_xX \text{ and } x\in X.
  \end{equation}
\end{remark}

\begin{example}\label{ex:retractions_in_riemannian_manifolds}
  Let $(X,g)$ be a Riemannian manifold. Then the exponential map
  $\exp_g:\mathcal{Z}_g\rightarrow X$ is defined in an open neighborhood
  $\mathcal{Z}_g$ of $Z_X\subset TX$; notice that when $g$ is a
  complete metric $\mathcal{Z}_g=TX$. Then, $R\CE \exp_g$ is a smooth
  retraction. Indeed, $\exp_g(0_x) = x$ and, for any $v_x\in T_xX$,
  \begin{equation*}
    \begin{split}
      \frac{d}{dt}\bigg|_{t=0} \exp_g|_x(tv_x) = T_{0_x}(\exp_g|_x)(v_x) = v_x. 
    \end{split}
  \end{equation*}
  The smoothness of $\exp_g$ and the computation of $T(\exp_g|_x)$ are
  in Proposition 5.19
  of~\cite{bo:lee-introduction_to_riemannian_manifolds}.
\end{example}

If $R:\jcZ\rightarrow X$ is a smooth retraction, we define its
associated \jdef{extended retraction}
$\ti{R}:\jcZ\rightarrow X\times X$ by $\ti{R}(v_x) \CE  (x,R(v_x))$.

Before we tackle an important property of $\ti{R}$, we recall a well
known fact about $T_{0_x}TX$ for $x\in X$. Fix $x\in X$ and define two
vector subspaces $H_x,V_x\subset T_{0_x}TX$ by
\begin{equation}\label{eq:factors_of_T_0_xTX_decomposition}
  H_x\CE T_{0_x}Z_X \stext{ and }
  V_x \CE  \left\{\frac{d}{dt}\bigg|_{t=0}t\, \delta x : \delta x\in T_xX\right\}
  = \ker(T_{0_x}\tau_X). 
\end{equation}
Notice that $T_xX \simeq H_x$ (via
$\gamma'(0)\in T_xX \mapsto \frac{d}{dt}\big|_{t=0}0_{\gamma(t)} \in
H_x$) and $T_xX\simeq V_x$ (via
$\delta x\in T_xX \mapsto \frac{d}{dt}\big|_{t=0} t\, \delta x$). It
is easy to see that $T_{0_x}TX = H_x\oplus V_x$.

Now we prove that the extended retraction $\ti{R}$ is a diffeomorphism
between open neighborhoods of $Z_X$ and $\Delta_X$ in $TX$ and
$X\times X$ respectively.

\begin{prop}\label{prop:extended_retraction}
  $\ti{R}$ maps the zero section $Z_X\subset TX$ bijectively onto the
  diagonal $\Delta_X\subset X\times X$. Furthermore, there are open
  subsets $V_R\subset TX$ and $W_R\subset X\times X$ containing $Z_X$
  and $\Delta_X$ respectively, such that
  $\ti{R}|_{V_R}^{W_R}:V_R\rightarrow W_R$ is a diffeomorphism.
\end{prop}

\begin{proof}
  The first assertion is immediate from point~\eqref{it:retraction_0}
  in Definition~\ref{def:retraction}. The idea for the second
  assertion is as follows: it is easy to see by direct computation
  using the Inverse Function Theorem (or invoking Proposition 2.1
  in~\cite{ar:barberoLinan_martinDeDiego-retratction_maps_a_seed_of_geometric_integrators}
  together with the fact that $\ti{R}$ is a discretization map as in
  Definition 2.2 of the same paper) that $\ti{R}$ is a local
  diffeomorphism in a neighborhood of $Z_X$. Then, the existence of
  the open subsets $V_R$ and $W_R$ follows using Theorem 1
  in~\cite{ar:cuell_patrick-skew_critical_problems}\footnote{The proof
    of this result refers to~\cite{bo:lang-differential_manifolds}. An
    alternative reference is to follow Exercise 14 on p. 56
    of~\cite{bo:Guillemin-Pollack-differential_topology}.}.
\end{proof}

From~\eqref{it:retraction_0} in Definition~\ref{def:retraction} we see
that $(\ti{R}|_{V_R}^{W_R})^{-1}(x,x) = 0_x$ for all $x\in X$. In
what follows, whenever we write $\ti{R}^{-1}$, we mean
$(\ti{R}|_{V_R}^{W_R})^{-1}$.


\subsection{Discrete connections from retractions and connections}
\label{sec:discrete_connections_from_retractions_and_connections}

Let $\phi:Q\rightarrow M$ be a principal $\SG$-bundle and $\CC$ be a
connection on $\phi$. Given a smooth retraction $R$ on $Q$, we want to
see how to ``discretize'' $\CC$ using $R$ in order to construct a
discrete connection $\DC$ on $\phi$ such that, ideally, it satisfies
$F_C(\DC)=\CC$. It is natural that, in order to do so, we may need to
impose some additional conditions to $R$.

The first step is the construction of a smooth retraction on $M$ using
$R$ and $\CC$.

As we know, the horizontal distribution of $\CC$, $Hor_\CC\subset TQ$,
is a $\SG$-invariant distribution on $Q$ that satisfies
$TQ=\VC\oplus Hor_\CC$, where $\VC\CE \ker(T\phi)$ is the \jdef{vertical
  bundle} of $\phi$. In particular, it is a regular submanifold
($Hor_\CC=\CC^{-1}(\{0\})$ and $0\in\jgg$ is a regular value of
$\CC$).

Assume that $R$ is $\SG$-equivariant for the actions $l^{TQ}$ and
$l^Q$. As $\phi\circ R$ is $\SG$-invariant, it induces a smooth map
$R_\SG:(TQ)/\SG\rightarrow M$ and we have the following commutative
diagram
\begin{equation*}
  \xymatrix{{(TQ)/\SG} \ar[r]^-{R_\SG} & {M} \ar@{<->}[d]^{id_{M}}\\
    {Hor_\CC/\SG} \ar@{^{(}->}[u]^i \ar[r]_-{\bar{R}}& {M}
  }
\end{equation*}
where $i$ is the natural inclusion and $\bar{R}$ is the restriction of
$R_\SG$ to $Hor_\CC/\SG$.

The connection $\CC$ allows us to define an isomorphism of vector
bundles over $M$ (see Section 2.4
in~\cite{bo:cendra_marsden_ratiu-lagrangian_reduction_by_stages}):
\begin{equation*}
  \alpha_\CC:(TQ)/\SG\rightarrow TM\oplus \ti{\jgg} \stext{ by }
  \alpha_\CC([v_q]) : =
  (T_q\phi(v_q),[(q,\CC(v_q))]),
\end{equation*}
where 
$\tilde\jgg\CE (Q\times \jgg)/\SG$ is the \jdef{adjoint bundle} (where
$\SG$ acts by the principal bundle action on $Q$ and the adjoint
representation on $\jgg$). Notice that
$\alpha_\CC(Hor_\CC/\SG) = TM\oplus 0_M$. Using this new isomorphism
we have the following commutative diagram in the category of
manifolds:
\begin{equation*}
  \xymatrix{{TM\oplus\ti{\jgg}} & {} & {(TQ)/\SG} \ar[ll]_-{\alpha_\CC}^-{\sim} \ar[r]^-{R_\SG} & {M} \ar@{<->}[d]^{id_{M}}\\
    {TM\oplus 0_{M}} \ar@{^{(}->}[u]^i & {} & {Hor_\CC/\SG} \ar[ll]^-{\alpha_\CC|_{Hor_\CC/\SG}}_-{\sim} \ar@{^{(}->}[u]^i \ar[r]_{\bar{R}}& {M} \\
    {TM} \ar[u]^{i_0}_{\sim} \ar@/_1pc/[rrru]_{\check{R}} & {} & {} & {}
  }
\end{equation*}
where we have defined
\begin{equation}\label{eq:reduced_retraction_abstract-def}
  \check{R}:TM\rightarrow M \stext{ by }
  \check{R} \CE  \bar{R}\circ (\alpha_\CC|_{Hor_\CC/\SG})^{-1} \circ i_0. 
\end{equation}
More explicitly, if $\delta m\in T_mM$ and $q\in Q|_m$,
\begin{equation}\label{eq:reduced_retraction_explicit-def}
  \check{R}(\delta m) \CE 
  \bar{R}((\alpha_\CC|_{Hor_\CC/\SG})^{-1}(\delta m,0)) =
  R_\SG([\HLc{}(q,\delta m)]) = \phi(R(\HLc{}(q,\delta m))),
\end{equation}
where $\HLc{}=F_{CH}(\CC)$ is the horizontal lift map associated to
the connection $\CC$.

\begin{lemma}
  If $R$ is a $\SG$-equivariant smooth retraction on $Q$, the map
  $\check{R}$ defined by~\eqref{eq:reduced_retraction_abstract-def} is
  a smooth retraction on $M$.
\end{lemma}

\begin{proof}
  Being a composition of smooth maps, $\check{R}$ is a smooth
  map. Also, for any $m\in M$, and $q\in Q|_m$,
  \begin{equation*}
    \check{R}(0_m) =  \phi(R((\HLc{}(q,0_m)))) =  \phi(R(0_q)) =\phi(q) = m.
  \end{equation*}
  Last, for $\delta m\in T_mM$ and $q\in Q|_m$,
  \begin{equation*}
    \begin{split}
      \frac{d}{dt}\bigg|_{t=0} \check{R}(t\, \delta m) =&
      \frac{d}{dt}\bigg|_{t=0} \phi(R(\HLc{}(q,t\, \delta m))) =
      \frac{d}{dt}\bigg|_{t=0} \phi(R(t\, \HLc{}(q,\delta m))) \\=&
      T_q\phi\bigg(\frac{d}{dt}\bigg|_{t=0} R(t\, \HLc{}(q,\delta
      m))\bigg) = T_q\phi(\HLc{}(q,\delta m)) = \delta m.
    \end{split}
  \end{equation*}
\end{proof}

As before, let $\CC\in\Sigma_C$ and $\HLc{}\CE F_{CH}(\CC)\in
\Sigma_H$. Define
\begin{equation}\label{HLd_from_CC_and_retraction-def-1}
  \HLd{}: \jcU'\rightarrow Q\times Q \stext{ by } \HLd{}(q,m) \CE 
  \ti{R}\big(\HLc{}(q,{(\ti{\check{R}})}^{-1}(\phi(q),m))\big),
\end{equation}
for
\begin{equation}
  \label{eq:HLd_from_CC_and_retraction-def-U}
  \jcU\CE (\phi\times\phi)^{-1}(W_{\check{R}}) \stext{ and }
  \jcU'\CE (id_Q\times \phi)(\jcU),
\end{equation}
where $W_{\check{R}}\subset M\times M$ is the open subset constructed
by Proposition~\ref{prop:extended_retraction} for $\check{R}$.

\begin{remark}\label{rem:HLd_from_R_is_well_defined}
  For $\HLd{}$ to be well defined, we need to check that, for
  $(q,m)\in\jcU'$,
  $(q,{\ti{\check{R}}}^{-1}(\phi(q),m)) \in \phi^*TM$, that is, we
  have to check that
  $\phi(q) = \tau_M({\ti{\check{R}}}^{-1}(\phi(q),m))$. This is true
  because for $\delta m\in TM$, such that
  $\delta m = {\ti{\check{R}}}^{-1}(\phi(q),m)$, we have
  $(\tau_M(\delta m),\check{R}(\delta m)) = \ti{\check{R}}(\delta m) =
  (\phi(q),m)$, so that
  $\tau_M({\ti{\check{R}}}^{-1}(\phi(q),m)) = \tau_M(\delta m)
  =\phi(q)$.
\end{remark}

\begin{prop}\label{prop:HLD_from_CC_and_retraction-1}
  Let $\phi:Q\rightarrow M$ be a principal $\SG$-bundle, $\CC$ be a
  connection on $\phi$ and $R$ be a $\SG$-equivariant smooth
  retraction on $Q$. Then,
  \begin{enumerate}
  \item \label{it:HLD_from_CC_and_retraction-1-D_type} the set $\jcU$
    defined in~\eqref{eq:HLd_from_CC_and_retraction-def-U} is of
    $D$-type, and
  \item \label{it:HLD_from_CC_and_retraction-1-HL} the map $\HLd{}$
    defined by~\eqref{HLd_from_CC_and_retraction-def-1} is in
    $\Sigma_H^d(\jcU)$.
  \end{enumerate}  
\end{prop}

\begin{proof}
  By Proposition~\ref{prop:extended_retraction},
  $\Delta_M\subset W_{\check{R}}$ and $W_{\check{R}}\subset M\times M$ is
  open. Thus, by the continuity of $\phi$, so are $\jcU$ and
  $\jcU'=(id_Q\times \phi)(\jcU) = (\phi\times
  id_M)^{-1}(W_{\check{R}})$. By the $\SG$-invariance of $\phi$,
  $\jcU$ is $\SG\times\SG$-invariant. As
  $\Delta_M\subset W_{\check{R}}$, we have $\Delta_Q\subset \jcU$ and
  we conclude that
  point~\eqref{it:HLD_from_CC_and_retraction-1-D_type} is
  true. Point~\eqref{it:HLD_from_CC_and_retraction-1-HL} essentially
  follows by unraveling the definitions in order to check that
  $\HLd{}$ is a discrete connection with domain $\jcU'$.
\end{proof}

\begin{theorem}\label{th:HLd_from_retraction_induces_HLc-1}
  Let $\phi:Q\rightarrow M$ be a principal $\SG$-bundle, $R$ be a
  $\SG$-equivariant smooth retraction on $Q$ and $\HLc{}$ a horizontal
  lift on $\phi$. Define $\HLd{}\in \Sigma_H^d(\jcU)$
  by~\eqref{HLd_from_CC_and_retraction-def-1}
  and~\eqref{eq:HLd_from_CC_and_retraction-def-U}. Then,
  $F_H(\HLd{})=\HLc{}$, that is, $\HLd{}$ integrates $\HLc{}$.
\end{theorem}

\begin{proof}
  Let $\widehat{\HLc{}}\CE F_H(\HLd{})$, we want to prove that
  $\widehat{\HLc{}}=\HLc{}$.  Fix $(q,\delta m) \in \phi^*TM$. Then,
  (see Proposition~\ref{prop:F_H-def}),
  \begin{equation*}
    \begin{split}
      \widehat{\HLc{}}(q,\delta m) =& (D_2\HLds{})(q,\phi(q))(\delta m) =
      (D_2(R(\HLc{}(q,\ti{\check{R}}^{-1}(\phi(q_0),m_1)))))
      \big|_{q_0=q,m_1=\phi(q)}(\delta m).
    \end{split}
  \end{equation*}
  Thus, if $\gamma:(-\epsilon,\epsilon)\rightarrow M$ is a smooth
  curve satisfying $\gamma(0)=\phi(q)$ and $\gamma'(0)=\delta m$, we
  have
  \begin{equation}\label{eq:HLd_from_retraction_induces_HLc-HLdc_explicit}
    \widehat{\HLc{}}(q,\delta m) =
    \frac{d}{dt}\bigg|_{t=0}
    R(\HLc{}(q,\ti{\check{R}}^{-1}(\phi(q),\gamma(t)))).
  \end{equation}

  We notice that, as
  $id_{W_{\check{R}}} =\ti{\check{R}}\circ \ti{\check{R}}^{-1}$, we
  have
  \begin{equation*}
    id_{T_{(\phi(q),\phi(q))}(M\times M)} =
    T_{0_{\phi(q)}}\ti{\check{R}} \circ T_{(\phi(q),\phi(q))} \ti{\check{R}}^{-1}.
  \end{equation*}
  Evaluating the previous identity at $(0,\delta m)$ we obtain
  \begin{equation*}
    \begin{split}
      (0,\delta m) =&
      T_{0_{\phi(q)}}\ti{\check{R}}(T_{(\phi(q),\phi(q))}
      \ti{\check{R}}^{-1}(0,\delta m)) \\=&
      T_{0_{\phi(q)}}\ti{\check{R}}\bigg(T_{(\phi(q),\phi(q))}
      \ti{\check{R}}^{-1}\bigg(\frac{d}{dt}\bigg|_{t=0}(\phi(q),\gamma(t))\bigg)\bigg)
      \\=&
      T_{0_{\phi(q)}}\ti{\check{R}}\bigg(\frac{d}{dt}\bigg|_{t=0}\ti{\check{R}}^{-1}(\phi(q),\gamma(t))\bigg)
    \end{split}
  \end{equation*}
  Writing $\delta m(t)\CE \ti{\check{R}}^{-1}(\phi(q),\gamma(t))$, we have
  \begin{equation*}
    (\phi(q),\gamma(t)) = \ti{\check{R}}(\delta m(t)) = (\tau_M(\delta
    m(t)),\check{R}(\delta m(t))),
  \end{equation*}
  so that $\delta m(t) \in T_{\phi(q)}M$ for all $t$. As a
  consequence,
  $\frac{d}{dt}\big|_{t=0}\ti{\check{R}}^{-1}(\phi(q),\gamma(t)) \in
  V_{\phi(q)}\subset T_{0_{\phi(q)}}T_{\phi(q)} M$
  (see~\eqref{eq:factors_of_T_0_xTX_decomposition}). Then,
  \begin{equation}\label{eq:HLd_from_retraction_induces_HLc-tiCheckRinv}
    \begin{split}
      (0,\delta m) =&
      T_{0_{\phi(q)}}\ti{\check{R}}\bigg(\frac{d}{dt}\bigg|_{t=0}\ti{\check{R}}^{-1}(\phi(q),\gamma(t))\bigg)
      \\=&
      \bigg(T_{0_{\phi(q)}}\tau_M\bigg(\frac{d}{dt}\bigg|_{t=0}\ti{\check{R}}^{-1}(\phi(q),\gamma(t))\bigg),
      T_{0_{\phi(q)}}\check{R}\bigg(\frac{d}{dt}\bigg|_{t=0}\ti{\check{R}}^{-1}(\phi(q),\gamma(t))\bigg)\bigg)
      \\=&
      \bigg(0,T_{0_{\phi(q)}}\check{R}_{\phi(q)}\bigg(\frac{d}{dt}\bigg|_{t=0}\ti{\check{R}}^{-1}(\phi(q),\gamma(t))\bigg)\bigg)
    \end{split}
  \end{equation}
  Now, by condition~\eqref{it:retraction_rigidity} in
  Definition~\ref{def:retraction},
  $T_{0_{\phi(q)}}\check{R}_{\phi(q)} = id_{T_{\phi(q)}}M$, but notice
  that this last identity assumes the identification
  $T_{0_{\phi(q)}}T_{\phi(q)}M$ with $T_{\phi(q)}M$. If we don't make
  this identification,
  from~\eqref{eq:HLd_from_retraction_induces_HLc-tiCheckRinv} we
  obtain
  \begin{equation*}
    \frac{d}{dt}\bigg|_{t=0} t\ \delta m =
    \frac{d}{dt}\bigg|_{t=0}\ti{\check{R}}^{-1}(\phi(q),\gamma(t)).
  \end{equation*}

  Now, we have
  \begin{equation*}
    \begin{split}
      \frac{d}{dt}\bigg|_{t=0}
      \HLc{}(q,(\ti{\check{R}}^{-1}(\phi(q),\gamma(t)))) =&
      T_{(q,\phi(q))}\HLc{} \bigg(
      \frac{d}{dt}\bigg|_{t=0}(q,(\ti{\check{R}}^{-1}(\phi(q),\gamma(t))))
      \bigg) \\=& T_{(q,\phi(q))}\HLc{}
      \bigg(0_q,\frac{d}{dt}\bigg|_{t=0}(\ti{\check{R}}^{-1}(\phi(q),\gamma(t)))
      \bigg) \\=& T_{(q,\phi(q))}\HLc{}
      \bigg(0_q,\frac{d}{dt}\bigg|_{t=0} t\ \delta m \bigg) =
      \frac{d}{dt}\bigg|_{t=0} \HLc{}(q,t\ \delta m) \\=&
      \frac{d}{dt}\bigg|_{t=0} t\ \HLc{}(q,\delta m) \in V_q \subset
      T_{0_q}T_qQ,
    \end{split}
  \end{equation*}
  where, in the last equality, we used the fact that $\HLc{}$ is
  (fiberwise) a linear map.

  Finally,
  using~\eqref{eq:HLd_from_retraction_induces_HLc-HLdc_explicit} and,
  later,~\eqref{eq:retraction_cond_2_explicit},
  \begin{equation*}
    \begin{split}
      \widehat{\HLc{}}(q,\delta m) =& \frac{d}{dt}\bigg|_{t=0}
      R(\HLc{}(q,\ti{\check{R}}^{-1}(\phi(q),\gamma(t)))) \\=&
      T_{0_q}R\bigg(
      \frac{d}{dt}\bigg|_{t=0}\HLc{}(q,(\ti{\check{R}}^{-1}(\phi(q),
      \gamma(t))))\bigg) = T_{0_q}R\bigg(\frac{d}{dt}\bigg|_{t=0}
      t\ \HLc{}(q,\delta m)\bigg) \\=& \frac{d}{dt}\bigg|_{t=0} R(t\
      \HLc{}(q,\delta m)) = \frac{d}{dt}\bigg|_{t=0} R_{q}(t\
      \HLc{}(q,\delta m))=\HLc{}(q,\delta m).
    \end{split}
  \end{equation*}
\end{proof}

\begin{corollary}\label{cor:DC_from_retraction_induces_CC}
  Let $\phi:Q\rightarrow M$ be a principal $\SG$-bundle, $R$ be a
  $\SG$-equivariant smooth retraction on $Q$, and $\CC\in
  \Sigma_C$. If $\HLd{}\in \Sigma_H^d(\jcU)$ is the one defined by
  Theorem~\ref{th:HLd_from_retraction_induces_HLc-1} using
  $\HLc{}\CE F_{CH}(\CC)$, then $\DC\CE (F_{CH}^d)^{-1}(\HLd{})$
  integrates $\CC$, that is, $F_C(\DC)=\CC$.
\end{corollary}

\begin{proof}
  By the commutativity of diagram~\eqref{eq:full_D_and_C_diagram} and
  Theorem~\ref{th:HLd_from_retraction_induces_HLc-1}, we have
  \begin{equation*}
    F_C(\DC)=F_C((F_{CH}^d)^{-1}(\HLd{})) = F_{CH}^{-1}(F_H(\HLd{})) =
    F_{CH}^{-1}(\HLc{}) = \CC.
  \end{equation*}
\end{proof}

Using Theorem~\ref{th:HLd_from_retraction_induces_HLc-1} or
Corollary~\ref{cor:DC_from_retraction_induces_CC}, the Integration
Problem has a solution, provided that a $\SG$-equivariant smooth
retraction on the principal $\SG$-bundle $\phi$ is available. The next
result provides a context where such retractions do exist.

\begin{theorem}\label{th:there_are_integrals_if_invariant_metrics}
  Let $\phi:Q\rightarrow M$ be a principal $\SG$-bundle and
  $\ip{}{}_Q$ be a $\SG$-invariant Riemannian metric on $Q$. Then, for
  any connection $\CC$ on $\phi$, there is, at least, one discrete
  connection $\DC$ on $\phi$ that integrates $\CC$, that is, it
  satisfies $F_C(\DC)=\CC$.
\end{theorem}

\begin{proof}
  As seen in Example~\ref{ex:retractions_in_riemannian_manifolds},
  $\ip{}{}_Q$ induces a smooth retraction $R$ on $Q$; explicitly,
  $R(\delta q) = \exp_{\ip{}{}_Q}(\delta q)$. As $\ip{}{}_Q$ is
  $\SG$-invariant, its geodesics are mapped to geodesics by the
  $\SG$-action $l^Q$; that is,
  $\exp_{\ip{}{}_Q}\circ l^{TQ}_g = l^Q_g\circ \exp_{\ip{}{}_Q}$, and
  it follows that $R$ is $\SG$-equivariant. If $\CC$ is any connection
  on $\phi$, the existence of a discrete connection $\DC$ on $\phi$
  such that $F_C(\DC)=\CC$ is guaranteed by
  Corollary~\ref{cor:DC_from_retraction_induces_CC}.
\end{proof}

Finally, the following result, based on a classical result on the
existence of invariant metrics on Lie groups, gives a sufficient
condition on the structure group of the principal bundle to ensure the
existence of invariant Riemannian metrics.

\begin{prop}\label{prop:existence_of_invariant_riemannian_metrics_on_PB}
  Let $\SG$ be a connected Lie group that is the Cartesian product of
  a compact Lie group and a vector space (with its additive group
  structure). Then, any principal $\SG$-bundle has a $\SG$-invariant
  Riemannian metric defined in its total space.
\end{prop}

\begin{proof}
  For $\SG$ as in the statement, let $\phi:Q\rightarrow M$ be a
  principal $\SG$-bundle. Pick a principal connection $\CC$ on $\phi$
  and a Riemannian metric $\ip{}{}^M$ on $M$. The horizontal
  distribution associated to $\CC$, $Hor_\CC\subset TQ$, is a
  $\SG$-invariant regular distribution satisfying
  $TQ=\VC\oplus Hor_\CC$, where $\VC\CE \ker(T\phi)$ is the
  \jdef{vertical bundle} of $\phi$. Let $Hor_\CC(q)$ be the fiber of
  $Hor_\CC$ over $q\in Q$; it is immediate that
  $T_q\phi|Hor_\CC(q):Hor_\CC(q)\rightarrow T_{\phi(q)}M$ is an
  isomorphism of vector spaces. For
  $\delta q, \delta q'\in Hor_\CC(q)$ define
  \begin{equation*}
    \ip{\delta q}{\delta q'}^{Hor_\CC}\CE 
    \ip{T_q\phi(\delta q)}{T_q\phi(\delta q')}^M.
  \end{equation*}
  It is easy to check that $\ip{}{}^{Hor_\CC}$ is a $\SG$-invariant
  Riemannian metric on the vector bundle $Hor_\CC\rightarrow Q$.

  Next, we want to find a $\SG$-invariant Riemannian metric on
  $\VC$. As vector bundles over $Q$, $Q\times\jgg \simeq \VC$ with the
  isomorphism given by $(q,\xi) \mapsto \xi_Q(q)$. The $\SG$-action
  $l^{TQ}$ induces a $\SG$-action on $\VC$:
  \begin{equation*}
    l^{TQ}_g(\xi_Q(q)) = (\Ad_g(\xi))_Q(l^Q_g(q)) \stext{ for all }
    q\in Q,\ \xi\in \jgg,\ g\in\SG. 
  \end{equation*}
  Using the previous isomorphism, this action induces the $\SG$-action
  on $Q\times \jgg$:
  \begin{equation*}
    l^{Q\times\jgg}_g(q,\xi)\CE (l^Q_g(q), \Ad_g(\xi)) \stext{ for all }
    q\in Q,\ \xi\in \jgg,\ g\in\SG.
  \end{equation*}
  Thus, in order to find a $\SG$-invariant metric on $\VC$, we need to
  find such a metric on $Q\times \jgg$ that is invariant for
  $l^{Q\times\jgg}$. This last task we can do if we find an
  $\Ad$-invariant inner product on $\SG$. It is well known that this
  last problem is equivalent to that of finding a bi-invariant metric
  on $\SG$. That those metrics exist is equivalent (see \cite[Lemma
  7.5]{ar:milnor-curvatures_of_left_invariant_metrics_on_lie_groups})
  to $\SG$ being connected and isomorphic to the Cartesian product of
  a compact Lie group with an additive vector group. As these
  conditions are satisfied in the statement, there is an
  $\Ad$-invariant metric in $\jgg$ which produces a $\SG$-invariant
  Riemannian metric $\ip{}{}^\VC$ on the vector bundle
  $\VC\rightarrow Q$.

  Finally, as $TQ=\VC\oplus Hor_\CC$, defining
  $\ip{}{}^Q\CE \ip{}{}^\VC\oplus \ip{}{}^{Hor_\CC}$ (making the two
  subbundles orthogonal) we have that $\ip{}{}^Q$ is a $\SG$-invariant
  Riemannian metric on $Q$, as required.
\end{proof}

\begin{remark}
  With a minor alteration in the proof, the result of
  Proposition~\ref{prop:existence_of_invariant_riemannian_metrics_on_PB}
  remains valid for structure groups that are compact but not
  necessarily connected. In this case, the existence result follows
  from~\cite[Proposition
  2.24]{bo:alexandrino_bettiol-lie_groups_and_geometric_aspects_of_isometric_actions}
  instead
  of~\cite{ar:milnor-curvatures_of_left_invariant_metrics_on_lie_groups}.
\end{remark}

\begin{corollary} \label{cor:existence_of_integrals_for_compact_G}
  If $\SG$ is a Lie group that is either
  \begin{enumerate}
  \item compact, or
  \item a Cartesian product of a compact and connected Lie group with
    a vector space,
  \end{enumerate}
  then the Integration Problem on any principal $\SG$-bundle has a
  solution, that is, for any $\CC\in\Sigma_C$, there is a $D$-type set
  $\jcU$ and $\DC\in\Sigma_C^d(\jcU)$ such that $F_C(\DC)=\CC$.
\end{corollary}

The previous results show that, in many circumstances, integrals of
principal connections do exist. Still, those integrals are far from
unique, as we see next.

\begin{example}
  Consider the principal bundle $p_1:\R^2\rightarrow \R$ with
  structure group $\R$ acting on $\R^2$ by translation on the second
  component. The $1$-form $\CC\CE dy$ defines a principal connection on
  $p_1$. On the other hand, let $f:\R^2\rightarrow \R$ be any smooth
  function; define $C:\R^2\rightarrow \R$ by
  $C(x_0,x_1)\CE (x_1-x_0)^2f(x_0,x_1)$. By
  Example~\ref{ex:trivial_budle:dc}
  $\DCp{f}((x_0,y_0),(x_1,y_1)) \CE  y_1-y_0 + C(x_0,x_1)$ is a discrete
  connection on $p_1$ with domain $\jcU\CE \R^2\times\R^2$. Then,
  \begin{equation*}
    \begin{split}
      F_C(\DCp{f})(\dot{x}\del_x+\dot{y}\del_y) =&
       D_2\DCp{f}((x,y),(x,y))(\dot{x}\del_x+\dot{y}\del_y) =
                                                   \frac{d}{dt}\bigg|_{t=0} (t\dot{y} + C(x,x+t\dot{x})) \\=&
                                                                                                              \dot{y} + \frac{d}{dt}\bigg|_{t=0} ((t\dot{x})^2f(x,x+t\dot{x})) = \dot{y},
    \end{split}
  \end{equation*}
  so that $F_C(\DCp{f}) = dy = \CC$. Thus, we have a (different)
  integral of $\CC$ for each different smooth function $f$.
\end{example}


\section{The abelian structure group case}
\label{sec:the_abelian_structure_group_case}

In this section, we specialize the Integration Problem to the case
where the structure group of the principal bundle is abelian. We see
first that the Integration Problem for arbitrary curvatures has a
solution if and only if the curvatures of the continuous and discrete
connections are related in a specific way. Then, we prove that when
that relation is realized, the Integration Problem for connections
whose curvatures satisfy that relationship admits a unique solution.

Throughout this section $\SG$ denotes a connected abelian Lie group and, in
keeping with the usual tradition, the product in $\SG$ is written with
the additive notation. Let us also  fix a principal $\SG$-bundle
$\phi:Q\rightarrow M$ until the end of this section. Unless stated
otherwise, both continuous and discrete connections are assumed to be
on~$\phi$.

\begin{proposition}
\label{prop:abelian:integration_exists}
  For any $\CC\in\Sigma_C$ there is a
  $D$-type set $\jcU$ and $\DC\in\Sigma_C^d(\jcU)$ such that
  $F_C(\DC)=\CC$.
\end{proposition}

\begin{proof}
  Any connected abelian Lie group can be expressed as a Cartesian
  product of a torus and a vector space --- see C.\,Procesi's
  book~\cite[\S4.2]{bo:procesi-lie_groups_an_approach_through_invariants_and_representations}
  --- and then
  Corollary~\ref{cor:existence_of_integrals_for_compact_G} yields that
  the Integration Problem has a solution.
\end{proof}

The following two lemmas show, both in the continuous and in the
discrete case, that the difference between connections on a principal
bundle give rise to connections on the trivial bundle.

\begin{lemma}
\label{l:abelian:descent:conti}
  Let $\CC,\CC':TQ\to\jgg$ be two connection forms on $\phi$. There
  exists a unique $\epsilon\in\Omega^1(M,\jgg)$ such that $\epsilon\circ
  T\phi = \CC-\CC'$. 
\end{lemma}

\begin{proof}
  Suppose that $q,q'\in Q$, $\delta q\in T_qQ$ and $\delta q'\in T_{q'}
  Q$ are such that $\phi(q) =\phi(q')$ and $T_q\phi(\delta q) =
  T_{q'}\phi (\delta q')$. Let $g\in G$ be such that $q'=l_g^Q(q)$. Consider
  $v\CE \delta q' - Tl^{Q}_g(\delta q)\in T_{q'}Q$ and observe
  that 
  \[
    T_{q'}\phi(v)
    = T_{q'}\phi(\delta q') - T(\phi \circ l^{Q}_g)_{q} (\delta q)
    = T_{q'}\phi(\delta q') - T_{q}\phi(\delta q) 
    = 0,
  \]
  that is, $v$ is vertical. Letting then $\xi\in\jgg$ be such that $v
  = \xi_Q(q')$ we see that
  \begin{equation}
  \label{eq:abelian:desc:v}
    \CC(v) 
    = \xi
    = \CC'(v).
  \end{equation}
  As the adjoint action is trivial in an abelian group
  $\CC(l^{TQ}_g(\delta q)) = \Ad_g(\CC(\delta q)) = \CC(\delta q)$.
  Remembering that $Tl^Q_g=l^{TQ}_g$ we obtain that
  \begin{equation}\label{eq:abelian:equinvariance}
    \CC(\delta q') - \CC(\delta q)
    = \CC \left(  \delta q' - l^{TQ}_g(\delta q) \right)
    = \CC(v).
  \end{equation}
  Now, this very equation~\eqref{eq:abelian:equinvariance} holds for
  $\CC'$, and putting this together with~\eqref{eq:abelian:desc:v} we
  see that
  \(
    \CC(\delta q') - \CC(\delta q)
    =
    \CC'(\delta q') - \CC'(\delta q).
  \)
  Rearranging, we obtain
  \begin{equation}
  \label{eq:abelian:desc:well-defined}
    \left(  \CC-\CC' \right)(\delta q)
    =
    \left( \CC-\CC' \right)(\delta q').
  \end{equation}

  Let us define for each open set $V$ in $M$ and section $\sigma:V\to
  Q$ of $\phi$ the function $\epsilon^\sigma:TV \to \jgg$ by
  $\epsilon^\sigma = \left( \CC-\CC' \right)\circ  T\sigma $.  Thanks
  to~\eqref{eq:abelian:desc:well-defined}, if $\sigma:V\to Q$ and
  $\sigma:V'\to Q$ are two such sections then the functions
  $\epsilon^\sigma$ and $\epsilon^{\sigma'}$ coincide in the
  intersections of their domains. This gives the desired function
  $\epsilon:TM\to\jgg$, and it is evidently linear at each point of $M$.
  The uniqueness is a consequence of the surjectivity of $T\phi$.
\end{proof}

\begin{lemma}\label{l:abelian:descent:disc}
  Let $\jcU\subset Q\times Q$ be of $D$-type and $\DC$, $\DC'$ two
  discrete connection forms with domain $\jcU$. There exists a unique
  smooth function $\zeta:\jcU''\CE (\phi\times\phi)(\jcU) \to G$ such that
  $\zeta\circ(\phi\times\phi)\mid_\jcU = \DC - \DC'$ and
  $\zeta(m,m)=0$ for every $m\in M$.
  
  Moreover, if $\DC$ and $\DC'$ have the same discrete curvature then
  the discrete connection form $\DCp \zeta$ on $M\times \SG $
  determined by $\zeta$ as in Example~\ref{ex:trivial_budle:dc} is
  flat.
\end{lemma}

\begin{proof}
  The difference $\difference\CE \DC-\DC':\jcU\to\SG$ satisfies
  that if $(q,q')\in\jcU$ and $g,g'\in\SG$ then 
  \begin{align*}
    \difference(l^Q_{g}(q),l^Q_{g'}(q'))
    &= \DC(l^Q_{g}(q),l^Q_{g'}(q')) - \DC'(l^Q_{g}(q),l^Q_{g'}(q')) 
    \\
    &= g' + \DC(q,q') - g - ( g' + \DC(q,q') - g)
    \\
    &= \difference(q,q').
  \end{align*}
  As a consequence of this invariance $\difference$ factors through
  the quotient $\jcU''$, which yields the desired $\zeta$.
  The uniqueness follows from the surjectivity of
  $(\phi\times\phi)\vert_{\jcU}$.

  Denote by $\DCp \zeta$ the discrete connection form on $M\times \SG$
  determined by $\zeta$ as in Example~\ref{ex:trivial_budle:dc}. If
  $(q_0,q_1,q_2)\in \jcU^{(3)}$ and $g_0,g_1,g_2\in G$, we write
  $m_i\CE\phi(q_i)$ for~$i=0,1,2$ and see that the curvature of
  $\DCp \zeta$ evaluated in
  $\left( (m_0,g_0),(m_1,g_1),(m_2,g_2) \right)$ is
  \begin{align*}
  \MoveEqLeft
    \DCp \zeta((m_0,g_0),(m_1,g_1)) - \DCp \zeta((m_0,g_0),(m_2,g_2)) 
    + \DCp \zeta((m_1,g_1),(m_2,g_2)) 
    \\
    & = g_1 + \zeta(m_0,m_1) - g_0
      - g_2 - \zeta(m_0,m_2) +  g_0
      + g_2 + \zeta(m_1,m_2) - g_1
    \\
    & = \alpha(q_0,q_1)
      - \alpha(q_0,q_2)
      + \alpha(q_1,q_2)
    = 0.
  \end{align*}
  This, of course, proves that $\DCp \zeta$ is flat.
\end{proof}

We next see that the uniqueness in the Integration Problem, whose
solution we already know exists thanks to
Proposition~\ref{prop:abelian:integration_exists}, can be controlled
through the curvature.

\begin{proposition}\label{prop:abelian:integration_unique}
  Let $\jcU\subset Q\times Q$ be of $D$-type and let $\DC$, $\DC'$ be
  two discrete connection forms on $\phi$ with domain $\jcU$. If $\DC$
  and $\DC'$ have the same discrete curvature and $F_C(\DC)=F_C(\DC')$
  then there exists a symmetric open subset of $\jcW$ of $D$-type on
  which $\DC$ and $\DC'$ agree.
\end{proposition}

\begin{proof}
  Write $\difference \CE  \DC - \DC'$ and, using
  Lemma~\ref{l:abelian:descent:disc}, let
  $\jcU''\CE (\phi\times\phi)(\jcU)$ and $\zeta:\jcU'' \to G$ such that
  $\zeta\circ(\phi\times\phi)\mid_\jcU = \DC - \DC'$. By virtue of
  Example~\ref{ex:trivial_budle:dc}, $\zeta$~determines a discrete
  connection form $\DCp{\zeta}$ on $p_1:M\times G\rightarrow M$
  defined on the open subset of $D$-type
  \begin{equation*}
    \jcU_0 \CE (p_1\times p_1)^{-1}(\jcU'') \subset
    (M\times G) \times (M\times G).
  \end{equation*}
  Given $q\in Q$ and $\delta q\in T_qQ$ we have
  \begin{equation*}
    \begin{split}
      D_2\zeta(\phi(q),&\phi(q))(T\phi(\delta q)) =
                         T_{(\phi(q),\phi(q))}\zeta (0,T\phi(\delta q))
      = T_{\left( q,q \right)} \alpha(0,\delta q)
        \\&=  T_{\left( q,q \right)}(\DC-\DC')(0,\delta q) =
      F_C(\DC)(\delta q) -F_C(\DC')(\delta q) =0.
    \end{split}
  \end{equation*}
  Thus, recalling Example~\ref{ex:trivial_budle:derivation},
  $\DCp{\zeta}$ integrates the trivial connection in $M\times G$.
  Thanks to the second part of Lemma~\ref{l:abelian:descent:disc}, the
  curvature of $\DCp{\zeta}$ is zero.  Consider now the discrete
  connection form $\DCp{0}$ on $M\times G$ defined by
  $\DCp{0}\left( (m_1,g_1),(m_2,g_2) \right) = g_2-g_1$ for any
  $m_1,m_2$ in $M$ and $g_1,g_2$ in $G$. It is straightforward to
  check that $\DCp{0}$ is a flat connection form that integrates the
  trivial continuous connection. Then we may replace $\jcU_0$ by its
  symmetrization and apply Theorem~\ref{th:solution_of_IP-flat_case}
  to find a symmetric open set $\jcU_1\subset \jcU_0$ of $D$-type such
  that the restrictions to $\jcU_1$ of $\DCp{\zeta}$ and $\DCp{0}$
  coincide. Let $\jcW\CE ((\phi,0)\times (\phi,0))^{-1}(\jcU_1)$. Then,
  if $(q_0,q_1)\in \jcW $, we have
  $\left((\phi(q_0),0),(\phi(q_1),0)\right)\in\jcU_1$ and, then,
  \begin{align*}
    \difference(q_0,q_1) 
    &
      = \zeta((\phi(q_0),\phi(q_1))
      = \DCp{\zeta}\left( (\phi(q_0),0), (\phi(q_1),0) \right)
    \\
    &
      = \DCp{0}\left( (\phi(q_0),0), (\phi(q_1),0) \right)
      = 0,
  \end{align*}
  which is tantamount to what we wanted to see. 
\end{proof}

In what remains of this section we see that fixing any discrete
curvature $\DB$ and a \emph{compatible} continuous curvature
$\RR(\DB)$, any connection $\CC$ with curvature $\RR(\DB)$ can be
integrated to a unique discrete connection $\DC$ whose curvature is
$\DB$.

By Proposition~\ref{prop:F_C_preserves_flatness}, we know that the
derivation process preserves the flatness of connections. The next
result, somehow, extends this result to the non-flat case, when $\SG$
is abelian.

\begin{proposition}\label{prop:abelian:samecurvatures}
  Let $\jcU\subset Q\times Q$ be of $D$-type and $\DC$, $\DC'$ be two
  discrete connection forms on $\phi$ with domain $\jcU$. If the
  discrete curvatures of $\DC$ and of $\DC'$ are equal then their
  derived connection forms $F_C(\DC)$ and $F_C(\DC')$ have the same
  curvature.
\end{proposition}

\begin{proof}
  Let $\difference\CE \DC-\DC':\jcU\to\SG$ and, using
  Lemma~\ref{l:abelian:descent:disc}, let
  $\zeta:\jcU''\CE (\phi\times\phi)(\jcU) \to G$ be the function such
  that $\zeta\circ(\phi\times\phi)\mid_\jcU = \difference$. We know
  from this same lemma that $\DCp{\zeta}$ is flat, so that, by
  Proposition~\ref{prop:F_C_preserves_flatness}, the associated
  connection $F_C(\DCp{\zeta})$ is also flat.
  Let $\eta$ be the form in $\Omega^1(M,\jgg)$ that determines the
  connection $F_C(\DCp \zeta)$ as in
  Example~\ref{ex:trivial_bundle:cc}: the fact that $F_C(\DCp \zeta)$
  is flat implies, thanks to
  Example~\ref{ex:trivial_bundle:flat:conti}, that $d\eta=0$.  Now,
  given $q\in Q$ and $\delta q\in T_qQ$,
  using Example~\ref{ex:trivial_budle:derivation},
  \begin{align*}
    \eta_{\phi(q)}(T\phi(\delta q)  ) 
    &
    = T_{\left( 
        \phi(q),\phi(q) \right)}\zeta\left( 0,T\phi(\delta q) 
      \right)
    \\
    &
    = T_{\left( q,q \right)}\alpha(0,\delta q)
    = F_C(\DC)(\delta q) - F_C(\DC') (\delta q)
  \end{align*}
  and, then, $\phi^\star\eta = F_C(\DC)-F_C(\DC')$. It follows that
  \begin{equation*}
    d\left( F_C(\DC)-F_C(\DC') \right) 
    = d\phi^\star\eta 
    = \phi^\star d\eta
    = 0,
  \end{equation*}
  and thus $dF_C(\DC)=dF_C(\DC')$, as desired.
\end{proof}

It follows from Proposition~\ref{prop:abelian:samecurvatures} that
given a discrete curvature form $\DB$ on an open set of $D$-type
$\jcU\subset Q\times Q$ ---this is, the curvature form of \emph{some}
discrete connection with domain $\jcU$---  the continuous curvature
form $\CB$ obtained by choosing a discrete connection $\DC$ with
curvature $\DB$ and computing the curvature of $F_C(\DC)$ is uniquely
determined. We write $\RR(\DB)\CE\CB$.

\begin{theorem}
  \label{thm:abelian:integrate_w_fixed_curvature}
  Let $\jcU\subset Q\times Q$ be of $D$-type and let $\DB$ be the
  discrete curvature of a discrete connection with domain $\jcU$. Given
  a continuous connection $\CC$ with curvature $\RR(\DB)$ there exists
  an open set of $D$-type $\jcU_0\subset\jcU$ and a discrete connection
  $\DC$ with domain $\jcU_0$ that integrates $\CC$ and has curvature
  $\DB$.
\end{theorem}

\begin{proof}
  Let $n\CE \dim \Hdr^1(M,\jgg)$ and
  $\eta_1,\ldots,\eta_n\in\Omega^1(M,\jgg)$ be closed differential forms
  whose classes in $\Hdr^1(M,\jgg)$ form a basis.

  Let $\DC'$ be a discrete connection with domain $\jcU$ and curvature
  $\DB$, and write $\CC'\CE F_C(\DC')$.  Following
  Lemma~\ref{l:abelian:descent:conti}, we let
  $\epsilon\in\Omega^1(M,\jgg)$ be such that
  $\epsilon\circ T \phi = \CC-\CC'$. Since $\CC$ and $\CC'$ have the
  same curvature we have that $d(\CC-\CC')=0$, and therefore the fact
  that $T\phi:TQ\to TM$ is surjective implies that $\epsilon$ is a
  closed differential form.  There exist then $c_1,\ldots,c_n\in\R$
  and $f:M\to\jgg$ such that
  \begin{equation*}
    \epsilon
    = df + c_1\eta_1+\cdots + c_n\eta_n, 
  \end{equation*}
  where the value of $df$ at each $\delta m\in T_mM$ is the element of
  $\jgg$ that identifies with $T_mf(\delta m)\in T_{f(m)}\jgg$ through
  the differential of the left multiplication by $f(m)$.  Let
  $i\in\left\{ 1,\dots,n \right\}$. Thanks to
  Example~\ref{ex:trivial_bundle:flat:conti} we may interpret
  $\eta_i\in\Omega^1(M,\jgg)$ as a flat connection $\CC_i$ on
  $p_1:M\times \SG\rightarrow M$, and use
  Theorem~\ref{th:solution_of_IP-flat_case} to find a subset $\jcU_i$
  of $(M\times \SG)\times (M\times \SG)$ of $D$-type and a flat
  discrete connection $\DCp{i}$ in $p_1:M\times\SG\rightarrow M$
  defined on $\jcU_i$ that integrates $\CC_i$.  If
  $\jcU_i''\CE (p_1\times p_1)(\jcU_i)$ then we have the local
  expression of $\DCp{i}$, $N_i:\jcU_i''\to \SG$
  (Example~\ref{ex:trivial_budle:dc}).  Evidently, $N_i(m,m)=0$ for
  all $m\in M$; by virtue of Example~\ref{ex:trivial_budle:derivation}
  we have that
  $T_{\left( m,m \right)}N_i(0,\delta m) = \eta_i(\delta m)$ whenever
  $m\in M$ and $\delta m \in T_mM$; moreover, using the flatness of
  $\DCp{i}$ it is straightforward to see that for any
  $m_0,m_1,m_2\in M$
  \begin{equation}\label{eq:local_discrete_flatness}
    0 = N_i(m_0,m_1) - N_i(m_0,m_2) +N_i(m_1,m_2).
  \end{equation} 
 
  Consider the function $\fint:M\times M$ defined by
  $\fint(m_1,m_2) \CE \exp_{\SG}\left( f(m_2)-f(m_1)
  \right)$. Set $\jcU_0'' \CE \bigcap_{i=1}^n\jcU_i''$ and
  define $\zeta:\jcU_0''\to \SG$ by
  $\zeta\CE \fint +\sum c_i N_i $.  We first observe that
  $\zeta(m,m)=0$ for every $m\in M$, and second that the map
  associated to $\zeta $ through the derivation functor is $\epsilon$:
  indeed, if $m\in M$ and $\delta m\in T_mM$ then
  \begin{align*}
    T_{\left( m,m \right)}\zeta(0,\delta m)
    & = T_{\left( m,m \right)}\fint(0,\delta m)
       + \sum c_i  T_{\left( m,m \right)}N_i(0,\delta m)
    \\
    & = T_mf(\delta m) + \sum c_i\eta_i(\delta m)
    = \epsilon(\delta m).
  \end{align*}
  Let $\jcU_0\CE (\phi\times\phi)^{-1}(\jcU''_0)\cap\jcU$ and
  $\difference\CE \zeta\circ(\phi\times\phi)\mid_{\jcU_0}$; if
  $m=\phi(q)$ and $\delta m=d\phi(\delta q)$ for $q\in Q$ and
  $\delta q \in T_qQ$ the calculation above means that
  \[
    T_{\left(q,q \right)}\alpha(0,\delta q)
    = T_{\left( m,m \right)}\zeta(0,\delta m)
    = \epsilon(\delta m )
    = \CC(\delta q) - \CC'(\delta q).
  \]
  Let now $\DC\CE \DC' + \alpha$. Evidently $\DC$ is a
  connection, for its value along the diagonal is null and
  \begin{align*}
    \DC(l^Q_{g}(q),l^Q_{g'}(q'))
    &
    = 
    \DC'(l^Q_{g}(q),l^Q_{g'}(q'))  + \alpha(l^Q_{g}(q),l^Q_{g'}(q'))
    \\
    &  
    = g' + \DC'(q,q') -g + \alpha(l^Q_{g}(q),l^Q_{g'}(q'))
    \\
    &
    = g' + \DC(q,q') -g
  \end{align*}
  for any $g,g'\in \SG$ and $(q,q')\in \jcU_0$. For any $q\in Q$ and
  $\delta q\in T_qQ$ we have
  \begin{align*}
    F_C(\DC)(  \delta q )
    &
      = T_{\left( q,q \right)}\DC'(0,\delta q) 
      + T_{\left(q,q \right)}\alpha(0,\delta q)
    \\
    &
      = \CC'(\delta q) + \CC(\delta q) - \CC'(\delta q)
      = \CC(\delta q).
  \end{align*}

  It is easy to see that the curvatures $\DB$ and $\DB'$ of $\DC$ and
  $\DC'$ are related as follows. For any
  $(q_0,q_1,q_2) \in \jcU_0^{(3)}$,
  \begin{equation*}
    \DB(q_0,q_1,q_2) =  \DB'(q_0,q_1,q_2) +
    \alpha(q_0,q_1) - \alpha(q_0,q_2) + \alpha(q_1,q_2).
  \end{equation*}
  If for any such $(q_0,q_1,q_2) \in \jcU_0^{(3)}$ we define
  $m_j\CE\phi(q_j)$ for $0\leq j \leq 2$, we have that
  \begin{align*}
    \MoveEqLeft 
    \alpha(q_0,q_1) - \alpha(q_0,q_2) + \alpha(q_1,q_2)
    \\
    &
      = 
      \zeta(m_0,m_1) - \zeta(m_0,m_2) + \zeta(m_1,m_2)
    \\
    & \!
    \begin{multlined}[t][.8\linewidth]
       = \fint (m_0,m_1) - \fint (m_0,m_2) + \fint (m_1,m_2)
       \\
       + \sum c_i \left( N_i(m_0,m_1) - N_i(m_0,m_2) +N_i(m_1,m_2)  \right).
    \end{multlined}
  \end{align*}
  The first three of these terms add up to zero --- this follows from
  the additivity of the exponential map in abelian groups --- and the
  other ones vanish in view of~\eqref{eq:local_discrete_flatness}.  It
  follows that $\DB=\DB'$, finishing the proof.
\end{proof}

\begin{corollary}
  Let $\jcU\subset Q\times Q$ be of $D$-type, $\DB$ be the discrete
  curvature of a discrete connection with domain $\jcU$ and $\CC$ be a
  continuous connection.  There exists an open subset
  $\openset\subset\jcU$ of $D$-type and a discrete connection $\DC$ with
  domain $\openset$ and curvature $\DB$ that integrates $\CC$ if and only
  if the curvature of $\CC$ is $\RR(\DB)$.

  In case these conditions hold and $\DC$ and $\DC'$ are two such
  discrete connections there exists a symmetric open subset of
  $\openset$ of $D$-type on which $\DC$ and $\DC'$ agree.
\end{corollary}

\begin{proof}
  The necessity is a consequence of the fact that $\RR(\DB)$ is well
  defined ---see Proposition~\ref{prop:abelian:samecurvatures}---, the
  sufficiency follows from
  Theorem~\ref{thm:abelian:integrate_w_fixed_curvature} above, and,
  finally, the uniqueness is guaranteed by
  Proposition~\ref{prop:abelian:integration_unique}.
\end{proof}


\appendix

\section{A basis of neighborhoods}
\label{sec:appendix}

Theorem~\ref{th:solution_of_IP-flat_case} ---one of our main
results--- relies on a particularization of the following claim made
in the proof of Proposition~2.\,1
in~\cite{ar:cabrera_marcut_salazar-on_local_integration_of_lie_brackets}:
\begin{quote}
  open sets of the form $U^{ (k)} =  G^{(k)} \cap U^ k$,
  where $U \subset G_\tau \cap G_\iota $ is an open neighborhood of $M$,
  form a basis of neighborhoods of $M$ in~$G^{(k)}$.
  \footnote{Here, $G$ is a Lie groupoid over a closed embedded submanifold
    $M\subset G$; the target, source and inversion maps $\tau$, $\sigma$ and
    $\iota$ have domains $G_\tau$, $G_\sigma$ and $G_\iota$; the superscript $k$
    denotes the $k$th power; and $G^{(k)} = \left\{ (g_1, \ldots, g_k) :
    \sigma(g_i) = \tau(g_{i+1}) \right\}$.}
\end{quote}
In response to our request, I. M\u{a}rcu\c{t} generously provided
guidelines for a more detailed explanation, which we further elaborate
below: the conclusions are formulated as Proposition~\ref{prop:basis}.
 
Let $X$ be a smooth manifold. Consider the manifold $P\CE X\times X$
and, for $j=1,2$, the fiber bundle $P^j$ over $X$ given by the map $p_j:P\to
X$ defined by $p_j(x_1,x_2) \CE x_j$. Our main proposition concerns
the fiber products $P_2\CE P^2 \prescript{}{p_2}{\times}_{p_1}P^1 $ and
$(\Delta_X)_2\CE \Delta_X \prescript{}{p_2}{\times}_{p_1}\Delta_X
\subset P_2$.

Let $g$ be a Riemannian metric on $X$, and
define $E_g^1,E_g^2 : TX\to P^j$ by 
$E_g^1(v_x) \CE (x,\exp_g(x))$ and
$E_g^2(v_x) \CE (\exp_g(x),x)$.
The following claims can be straightforwardly verified:

\begin{enumerate}
  \item 
  There exists an open neighborhood $\Epsilon_g\subset TX$ of the zero
  section $0_X$ such that if $LP^j \CE E^j_g(\Epsilon_g)$ then
  $E^j_g:\Epsilon_g\to LP^j$ is a fiber-preserving diffeomorphism whenever
  $j=1,2$.

  \item 
  The map $E_g:TX\times TX\to P^2\times P^1$ defined by $E_g(v_x,w_y)
  \CE (E^2_g(v_x),E^1_g(w_y))$ restricts to a diffeomorphism from $\left(
  \Epsilon_g\times \Epsilon_g \right)\cap(TX\oplus TX)$ onto $\left( LP^2\times
  LP^1 \right)\cap P_2$, that are open neighborhoods of $0_X\oplus 0_X$ and of
  $(\Delta_X)_2$, respectively. 
\end{enumerate}

\begin{lemma}\label{lemma:opensofzero}
  Let $\vbundle\to X $ be a vector bundle with a bundle metric and
  write
  $\vbundle_\rho =\left\{ v_x\in \vbundle_x : \norm{v_x} < \rho(x)
    \text{ and }x\in X \right\}$ for each real-valued function $\rho$
  on $X$.  The family
  $\left\{ \vbundle_\rho : \text{$\rho:X\to\R_{>0}$ continuous}
  \right\}$ is a basis of open neighborhoods of the zero section
  of~$\vbundle$.
\end{lemma}

\begin{proof}
  This can be proved using Theorem 2.1.3 on p.\,56 of
  \cite{bo:jost-riemannian_geometry_and_geometric_analysis-7e}.
\end{proof}

\begin{proposition}\label{prop:basis}
  Let $W\subset P_2$ be an open neighborhood of $(\Delta_X)_2$. There
  exists an open neighborhood $U\subset P$ of $\Delta_X$ such that
  $U^{(2)}\CE (U\times U)\cap P_2$ is contained in $W$ and contains 
  $(\Delta_X)_2$.
\end{proposition}

\begin{proof}
  We may assume that $W\subset LP^2\times LP^1$ and define $W'\CE
  E_g^{-1}(W)$, which is contained in $\left( \Epsilon_g\times \Epsilon_g
  \right)\cap(TX\oplus TX)$ and contains $0_X\times 0_X$.

  The Riemannian metric $g$ on $X$ provides a bundle metric
  $\langle\cdot,\cdot\rangle$ on $TX$ and, then,  the law
  \(
    \langle(v_x,w_x),(v_x',w_x')\rangle 
    \CE 
    \langle v_x,v'_x\rangle
    + \langle w_x,w'_x\rangle
  \)
  endows $TX\oplus TX\to X$ with a bundle metric.  By
  Lemma~\ref{lemma:opensofzero} there exists $\rho:X\to\R_{>0}$ such that
  $(TX\oplus TX)_\rho\subset W'$.  Let $\rho' \CE
  {\tfrac{1}{\sqrt{2}}\rho}$ and consider $(TX)_{\rho'}$. By definition,
  $(TX)_{\rho'}\subset TX$ is an open neighborhood of $0_X$; on the other hand,
  a short computation shows that
  \begin{equation}
  \label{eq:claim:sum}
    \text{
    if $v_x,v_x'\in(TX)_{\rho'}$ then
    $(v_x,v_x')\in(TX\oplus TX)_\rho\subset W'$.
    }
  \end{equation}

  Define $U^j \CE E^j_g( (TX)_{\rho'} ) \subset P^j$ for $j=1,2$ and
  $U\CE U^2\cap U^1$.  It is immediate that $U$ is an open neighborhood of
  $\Delta_X$ in $P$, and one can see using \eqref{eq:claim:sum} that
  $U^{(2)}\CE \left( U\times U \right)\cap P_2$, that evidently contains
  $(\Delta_X)_2$, is contained in $W$.
\end{proof}


\printbibliography


\end{document}